\documentclass[leqno,12pt]{article}

\usepackage{amssymb}
\usepackage{euscript}
\usepackage[dvips]{graphicx}
\usepackage{flafter}
\evensidemargin\oddsidemargin
\begin{document}

\baselineskip=18pt
\setcounter{page}{1}

\renewcommand{\theequation}{\thesection.\arabic{equation}}
\newtheorem{theorem}{Theorem}[section]
\newtheorem{lemma}[theorem]{Lemma}
\newtheorem{proposition}[theorem]{Proposition}
\newtheorem{corollary}[theorem]{Corollary}
\newtheorem{remark}[theorem]{Remark}
\newtheorem{fact}[theorem]{Fact}
\newtheorem{problem}[theorem]{Problem}
\newtheorem{counterexample}[theorem]{Counterexample}

%%%%%%%%%%%%%%%%%%%%%%%%%%% Equation numberings
\newcommand{\eqnsection}{
\renewcommand{\theequation}{\thesection.\arabic{equation}}
    \makeatletter
    \csname  @addtoreset\endcsname{equation}{section}
    \makeatother}
\eqnsection
%%%%%%%%%%%%%%%%%%%%%%%%%%%

%%%%%%%%%%%%%% Bbb characters
%%%%%%%%%%%%%% Real numbers
\def\r{{\mathbb R}}
%%%%%%%%%%%%%% Expectation
\def\e{{\mathbb E}}
%%%%%%%%%%%%%% Probability
\def\p{{\mathbb P}}
%%%%%%%%%%%%%% Integers
\def\z{{\mathbb Z}}
%%%%%%%%%%%%%% Positive Integers
\def\n{{\mathbb N}}
%%%%%%%%%%%%%% Favourite set
\def\V{{\mathbb V}}
%%%%%%%%%%%%%%

%%%%%%%%%%%%%%%% Special symbols
%%%%%%%%%%%%%% Exponential
\def\ee{\mathrm{e}}
%%%%%%%%%%%%%% Differentiation
\def\d{\, \mathrm{d}}
%%%%%%%%%%%%%%
\def\esss{\mathop{\mathrm{esssup}}}

%%%%%%%%%%%%%%%%%%%%%%%%%%%%%%%%%%%%%%%%%%%%
%%%%%%%%%%%%%%%%%%%%%%%%%%%%%%%%%%%%%%%%%%%%
%%%%%%%%%%%%%%%%%%%%%%%%%%%%%%%%%%%%%%%%%%%%

%%%%%%%%%%%%%% Beginning of the text

\vglue80pt

\centerline{\Large\bf Upper limits of Sinai's walk in random
scenery}

\bigskip
\bigskip
\bigskip

\centerline{by}

\bigskip
\bigskip

\centerline{Olivier Zindy}

\medskip

\centerline{\it Universit\'e Paris VI}

\bigskip
\bigskip
\bigskip
\bigskip

{\leftskip=2truecm
\rightskip=2truecm
\baselineskip=15pt
\small

\noindent{\slshape\bfseries Summary.} We consider Sinai's walk in
i.i.d.\ random scenery and focus our attention on a conjecture of
R\'ev\'esz concerning the upper limits of Sinai's walk in
random scenery when the scenery is bounded from above. A close study
of the competition between the concentration property for Sinai's
walk and negative values for the scenery enables us to prove that
the conjecture is true if the scenery has ``thin" negative tails and
is false otherwise.

\bigskip


\bigskip


} %%%%%% End of narrower

\bigskip
\bigskip

\section{Introduction}
   \label{rwrers:s:intro}

\subsection{Random walk in random environment}
\label{rwrers:Sinaiswalkintro}

Problems involving random environments arise in different domains of
physics and biology.  Originally, one-dimensional random walk in
random environment appeared as a simple model for DNA transcription.
In the following, we consider the elementary model of
one-dimensional Random Walk in Random Environment (RWRE), defined as
follows. Let $\omega:=(\omega_i, \, i \in \z)$ be a family of
independent and identically distributed (i.i.d.) random variables
defined on $\Omega,$ which stands for the random environment. Denote
by $P$ the distribution of $\omega$ and by $E$ the corresponding
expectation.

Conditioning on $\omega$ (i.e., choosing an environment), we define
the RWRE $(X_n, \, n \ge 0)$ as a nearest-neighbor random walk on
$\z$ with transition probabilities given by $\omega$: $(X_n, \, n
\ge 0)$ the Markov chain satisfying $X_0=0$ and for $n \ge 0,$
$$
P_\omega \{ X_{n+1} = x+1 \, | \, X_n =x\} = \omega_x = 1- P_\omega
\{ X_{n+1} = x-1 \, | \, X_n =x\}.
$$

\noindent We denote by  $P_{\omega}$ the law of $(X_n, \, n \ge 0),$
by $E_{\omega}$ the corresponding expectation, and by $\p$ the joint
law of $(\omega,(X_n)_{n \ge 0})$. We refer to Zeitouni
\cite{zeitouni} for an overview of random walks in random
environment.

Throughout the paper, we make the following assumptions on $\omega$:
\begin{eqnarray}
   \exists \, \delta \in (0,1/2): \qquad
    P \{ \delta\le \omega_0 \le 1-\delta\} &=&1,
    \label{rwrers:ellipticity}
    \\
    E [\, \log ({1-\omega_0\over \omega_0})\, ]
 &=& 0,
    \label{rwrers:E(log)=0}
    \\
    \sigma^2 := \hbox{\rm Var} [\, \log ({1-\omega_0\over
    \omega_0})\, ]
 &>&0.
    \label{rwrers:Var>0}
\end{eqnarray}

\noindent Assumption (\ref{rwrers:ellipticity}) implies that $|\log
({1-\omega_0\over \omega_0})|$ is, $P$-a.s., bounded by the constant
$L := \log ({1-\delta\over \delta})$. It is a technical assumption,
which can be replaced by an exponential moment of $\log
({1-\omega_0\over \omega_0})$. According to a recurrence-transience
result due to Solomon \cite{solomon}, assumption
(\ref{rwrers:E(log)=0}) ensures that $(X_n)_{n \ge 0}$ is
$\p$-almost surely recurrent, i.e., the random walk hits any site
infinitely often. Assumption (\ref{rwrers:Var>0}) excludes the case
of deterministic environment, which corresponds to the homogeneous
symmetric random walk.

Under assumptions (\ref{rwrers:ellipticity})--(\ref{rwrers:Var>0}),
the RWRE is referred to as Sinai's walk. Sinai \cite{sinai} proves
that $X_n / (\log n)^2$ converges in law, under $\p$, toward a
non-degenerate random variable, whose distribution is explicitly
computed by Kesten \cite{kesten86} and Golosov \cite{golosov2}. This
result contrasts with the usual central limit theorem which gives
the convergence in law of $X_n / \sqrt{n}.$

Let
\begin{eqnarray}
    L(n,x)
 &:=& \# \left\{ 0\le i\le n: \; X_i =x\right\}, \qquad n \ge 0, \;
    x\in \z,
    \label{rwrers:xi}
    \\
    L(n,A)
 &:=& \sum_{x \in A} L(n,x), \qquad n \ge 0, \;
    A \subset \z.
    \label{rwrers:L}
    \nonumber
\end{eqnarray}

\noindent In words, the quantity $L(n,A)$ measures the number of
visits to the set $A$ by the walk in the first $n$ steps.

The maximum of local time is studied by R\'ev\'esz (\cite{revesz05},
p.~$337$) and Shi \cite{shi98}: under assumptions
(\ref{rwrers:ellipticity})--(\ref{rwrers:Var>0}), there exists
$c_{0}>0$ such that
$$
\limsup_{n \to \infty} {\max_{x \in \z} L(n,x) \over n} \ge c_{0},
\qquad
 \p \textrm{-a.s.}
$$

\noindent It means that the walk spends, infinitely often,  a
positive part of its life on a single site. The liminf behavior is analyzed by Dembo,
Gantert, Peres and Shi \cite{dembo-gantert-peres-shi}, who prove that
$$
\liminf_{n \to \infty} {\max_{x \in \z} L(n,x) \over n/ \log \log \log n} = c'_{0},
\qquad
 \p \textrm{-a.s.,}
$$

\noindent for some $c'_{0} \in (0,\infty).$  A concentration
property is obtained by Theorem $1.3$ of Andreoletti
\cite{andreoletti}, which says that, under assumptions
(\ref{rwrers:ellipticity})--(\ref{rwrers:Var>0}) and for any
$0<\beta<1$, there exists $\ell(\beta)>0$ such that
\begin{equation}
\limsup_{n \to \infty} \, {\sup_{x \in \z}
L\left(n,\left[x-\ell(\beta),x+\ell(\beta)\right]\right) \over n } \ge
\beta, \qquad
 \p \textrm{-a.s.}
 \label{rwrers:andreo}
\end{equation}

\noindent In words, for any $\beta$ close to $1$, it is possible to
find a length $\ell(\beta)$ such that,  $\p$-almost surely, the
particle spends, infinitely often, more than a $\beta$-fraction of
its life in an interval of length $2 \ell(\beta).$

\medskip

\subsection{Random walk in random scenery}
\label{rwrers:rwrsintro}

Random Walk in Random Scenery (RWRS) is a simple model of diffusion
in disordered media, with long-range correlations. It is a class of
stationary random processes exhibiting rich behavior. It can be
described as follows: given a Markov chain on a state space, there
may be a random field indexed by the state space, called a random
scenery. As the random walk moves on this state space, he observes
the scenery at his location. For a survey of recent results about
RWRS, we refer to den Hollander and Steif \cite{denhollander-steif},
and to Asselah and Castell \cite{asselah-castell} for large
deviations results in dimension $d \ge 5.$

Let us now define the model of one-dimensional RWRS: consider
$S=(S_n, \, n\ge 0)$ a random walk on $\z$ and $\xi:=(\xi(x), \, x
\in \z)=(\xi_x, \, x \in \z)$, a family of i.i.d.\ random variables
defined on a probability space $\Xi$. We refer to $\xi$ as the
random scenery and denote by $Q$ its law. Then, define the process
 $(Y_n, \, n\ge 0)$ by
 $$
Y_n:=\sum_{i=0}^{n} \xi(S_i),
 $$

\noindent called RWRS or the Kesten-Spitzer Random Walk in Random
Scenery. An interpretation is the following: if a random walker has
to pay $\xi_x$ each time he visits $x$, then $Y_n$ stands for the
total amount he has paid in the time interval $[0,n]$.

The model is introduced and studied by Kesten and Spitzer
\cite{kesten-spitzer} in dimensions $d=1$ and $d \ge 3$. They prove
in dimension $d=1$ that, when $S$ and $\xi$ belong to the domains of
attraction of stable laws of indice $\alpha$ and $\beta$
respectively, then there exists $\delta,$ depending on $\alpha$ and
$\beta,$ such that $n^{-\delta} \, Y_{\lfloor nt \rfloor}$ converges
weakly. In the simple case where $\alpha=\beta=2,$ they show that
$$
\left(n^{-3/4} \, Y_{\lfloor nt \rfloor}; \, 0 \le t \le 1 \right)
\stackrel{\rm{law}}{\longrightarrow} \left( \Lambda(t) ; \, 0 \le t
\le 1 \right),
$$

\noindent where $``\stackrel{\rm{law}}{\longrightarrow}"$ stands for
weak convergence in law (in some functional space; for example in
the space of bounded functions on $[0,1]$ endowed with the uniform
topology). The process $(\Lambda(t), \, t \ge 0),$ called Brownian
motion in Brownian scenery, is defined by $\Lambda(0)=0$ and
$\Lambda(t):=\int_{\r}\ell(t,x) \d W(x)$ for $t>0,$ where $(W(x); \,
x \in \r)$ denotes a two-sided Brownian motion and $(\ell(t,x), \, t
\ge 0, x \in \r)$ denotes the jointly continuous version of the
local time process of a Brownian motion $(B(t), \, t \ge0),$
independent of $(W(x); \, x \in \r).$

Independently, Borodin analyzes the case of one-dimensional
nearest-neighbor random walk in random scenery, see \cite{borodinA}
and \cite{borodinB}. Bolthausen \cite{bolthausen} studies the case
$d = 2$. He proves that, if $S$ is a recurrent random walk  and
$\xi_0$ has zero expectation and finite variance, then $(n \log
n)^{-1/2} \, Y_{\lfloor nt \rfloor}$ satisfies a functional central
limit theorem.

\medskip

\subsection{Random environment and random scenery}
\label{rwrers:rersintro}

In this paper, we consider Sinai's Walk in Random Scenery. Problems
combining random environment and random scenery have been examined
for more general models. Replacing $\z$ by a more general countable
state space, Lyons and Schramm \cite{lyons-schramm} exhibit, under
certain conditions, a stationary measure for Random Walks in a
Random Environment with Random Scenery (RWRERS) from the viewpoint
of the random walker. H\"aggstr\"om \cite{haggstrom}, H\"aggstr\"om
and Peres \cite{haggstrom-peres} treat the case where the scenery
arises from percolation on a graph. In this particular case, the
scenery determines the random environment of the associated RWRE,
which is used by the authors to obtain information about the
scenery.

Let us first describe the model of Sinai's walk in random scenery.
We consider Sinai's walk $(X_n, \, n\ge0)$ under assumptions
(\ref{rwrers:ellipticity})--(\ref{rwrers:Var>0}), and recall that
the environment $\omega$ is defined on $(\Omega, P).$ For the
scenery, we consider a family of i.i.d.\ random variables
$\xi:=(\xi(x), \, x \in \z)=(\xi_x, \, x \in \z)$, defined on
$(\Xi,Q),$ independent of $\omega$ and $(X_n, \, n\ge0)$. To
translate independence between $\omega$ and $\xi$, we consider the
probability space $(\Omega \times \Xi,P\otimes Q)$, on which we
define $(\omega,\xi)$. Moreover, we denote by $\p \otimes Q$ the law
of $(\omega,(X_n)_{n \ge0},\xi)$. Then we define as Sinai's walk in
random scenery the process
 $(Z_n, \, n\ge 0)$:
 $$
Z_n:=\sum_{i=0}^{n} \xi(X_i).
 $$

\noindent Observe that $Z_n$ can be written using local time
notation:
\begin{equation}\label{rwrers:Zdecompo}
Z_n=\sum_{x \in \z} \xi(x) \, L(n,x), \qquad n \ge 0,
\end{equation}

\noindent where $L(n,x)$ stands for the local time of the random
walk at site $x$ until time $n,$ see (\ref{rwrers:xi}).

We are interested in the upper limit of $Z_n$ in the case where
$a:=\mathrm{ess } \sup \xi_0$ is finite. We consider the
concentration property of order $\beta$ for Sinai's walk with
$\beta$ close to $1$ (see (\ref{rwrers:andreo})), which enables us
to formulate the conjecture of R\'ev\'esz (\cite{revesz05},
p.~$353$): does the assumption that $a:=\mathrm{ess } \sup \xi_0$ is
finite imply that, $\p \otimes Q$-almost surely,
$$
 \limsup_{n \rightarrow \infty}{Z_n \over n}=a \, ?
$$

It turns out that the conjecture holds only under some additional
assumptions on the distribution of the random scenery. It is
interesting to note that this conjecture follows immediately from
the result of Andreoletti \cite{andreoletti} mentioned earlier if
$\mathrm{ess} \inf \xi_0$ is larger than $-\infty.$  In the general
case, a close study of the competition between the concentration
property for Sinai's walk and negative values for the scenery
enables us to obtain the following theorem, which gives a solution
to this problem, depending on the tail decay of $\xi^-_0:= \max
\{-\xi_0, 0 \}.$

\medskip

\begin{theorem}
 \label{rwrers:t:main}
 Assume (\ref{rwrers:ellipticity})--(\ref{rwrers:Var>0}) and $a:=\mathrm{ess } \sup \xi_0 < \infty.$
 \begin{itemize}
   \item[(i)]  If \,  $Q\{\xi^-_0 > \lambda\} \le {1 \over (\log \lambda )^{2+ \varepsilon}},$ for some
$\varepsilon>0$ and all large $\lambda$, then
\begin{eqnarray*}
    \p \otimes Q \left\{\limsup_{n \rightarrow \infty}{Z_n \over n}=a\right\}
    =1.
\end{eqnarray*}
     \item[(ii)]  If  \,  $Q\{\xi^-_0 > \lambda\} \ge {1 \over (\log \lambda)^{2-\varepsilon}},$ for some
$\varepsilon>0$ and all large $\lambda,$ then
\begin{eqnarray*}
  \p \otimes Q \left\{  \lim_{n \rightarrow \infty}{Z_n \over n}=- \infty \right\} =1.
\end{eqnarray*}
 \end{itemize}
\end{theorem}

\medskip

It is possible to give more precision in the case $(ii),$ see Remark
\ref{rwrers:r:negativpart}.  On the other hand, the case
$\varepsilon=0$ is still open.

We mention that, under
(\ref{rwrers:ellipticity})--(\ref{rwrers:Var>0}) and $a:=\mathrm{ess
} \sup \xi_0 < \infty,$ it is possible to prove that $\p \otimes Q
\{\limsup_{n \rightarrow \infty}{Z_n \over n}=c\},$ satisfies a
$0-1$ law, for any $c \in [-\infty,\infty].$ The proof follows the
lines of \cite{gantert-shi}, except that we need an additional
argument saying that modifying a finite number of random variables
in the scenery does not change the behavior of $\limsup_{n
\rightarrow \infty}{Z_n \over n}.$ The latter can be done by means
of Theorem $1$ (see Example $1$) in \cite{hu-shi98b}, which implies,
for any $x \in \z,$ that $L(n,x)=o(n),$ $n \to \infty,$ $\p$-almost
surely.

In general, we do not know whether $\limsup_{n \rightarrow
\infty}{Z_n \over n} \in \{-\infty,a\},$ $\p$-almost surely.

The paper is organized as follows: in Section
\ref{rwrers:s:preliminaries}, we present some key results for the
environment and for Sinai's walk when the environment is fixed
(i.e., quenched results). In Section \ref{rwrers:s:stepA}, we define
precisely the notion of ``good" environment-scenery and prove
Theorem \ref{rwrers:t:main} by accepting two intermediate
propositions. The first one, proved in Section
\ref{rwrers:proofpart}, is devoted to the study of the RWRE within
the ``good" environment-scenery. The second one, proved in Section
\ref{rwrers:proofenv}, does not concern the RWRE, but only the
environment-scenery. We show that, $P \otimes Q$-almost surely,
$(\omega,\xi)$ is a ``good" environment-scenery.

In the following, we use $c_i$ ($1\le i\le 33$) to denote finite and
positive constants.

\medskip

\section{Preliminaries}
\label{rwrers:s:preliminaries}

In this section, we collect some basic properties of random walk in
random environment that will be useful in the forthcoming sections.

\medskip

\subsection{About the environment}
\label{rwrers:subs:env}

In the study of one-dimensional RWRE, an important role is played by
a function of the environment $\omega$, called the potential. This
process, noted $V= (V(x), \; x\in \z)$, is defined on $(\Omega,P)$
by:
\begin{eqnarray}
\label{rwrers:potential}
 V(x) :=\left\{\begin{array}{lll} \sum_{i=1}^x \log
({1-\omega_i \over \omega_i}) & {\rm if} \ x \ge 1,
\\
 0 &  {\rm if} \ x=0,
\\
-\sum_{i=x+1}^0 \log ({1-\omega_i \over \omega_i}) &{\rm if} \ x\le
-1.
\end{array}
\right. \end{eqnarray}

\noindent By (\ref{rwrers:ellipticity}), we observe that
$|V(x)-V(x-1)| \le L$ for any $x\in \z.$ Moreover, we define $P_z
\{\cdot\} := P \{ \, \cdot \, | \, V(0)=z\}$, for any $z\in \r$;
thus $P=P_0$. (Strictly speaking, we should be
 working in a canonical space for $V$, with $P_z$ defined as the image measure
  of $P$ under translation.)

Let us define, for any Borel set $A\subset \r$,
\begin{eqnarray*}
    \nu^+(A) &:=& \min\left\{ n\ge 0: \, V(n) \in A\right\}.
    \label{rwrers:d+A}
\end{eqnarray*}

\noindent We recall the following result, whose proof is given by a
simple martingale argument.
\medskip

\begin{lemma}
 \label{rwrers:martingale1}
For any $x<y<z,$ we have
\begin{eqnarray*}
 \frac{y-x}{z-x+L}  \le P_y\{ \nu^+([z,\infty))< \nu^+((-\infty,x])\}
 &\le& \frac{y-x+L}{z-x}.
\end{eqnarray*}
\end{lemma}

\medskip

\noindent {\it Proof.} Since (\ref{rwrers:ellipticity}) and
(\ref{rwrers:E(log)=0}) imply that $(V(n); \, n \ge 0)$ is a
martingale with bounded jumps, we apply the optional stopping
theorem (\cite{durrett}, p.~$270$) at $\nu^+([z,\infty)) \wedge
\nu^+((-\infty,x])$ to get
\begin{eqnarray*}
y=E_y[X_0]&=&E_y[X_{\nu^+([z,\infty))} \, ; \, \nu^+([z,\infty))<
\nu^+((-\infty,x])]
\\
&&+E_y[X_{\nu^+((-\infty,x])} \, ; \, \nu^+([z,\infty))>
\nu^+((-\infty,x])].
\end{eqnarray*}

\noindent Since $X_{\nu^+([z,\infty))}\in[z,z+L]$ and
$X_{\nu^+((-\infty,x])} \in [x-L,x]$ by ellipticity, we obtain
$$
y \ge z P_y\{ \nu^+([z,\infty))< \nu^+((-\infty,x])\}+(x-L)(1-P_y\{
\nu^+([z,\infty))< \nu^+((-\infty,x])\}),
$$
 \noindent which yields the right inequality. The left
 inequality is obtained by similar arguments.
\hfill$\Box$

\bigskip

Moreover, we recall a result of Hirsch \cite{hirsch}, which, under
assumptions (\ref{rwrers:ellipticity})--(\ref{rwrers:Var>0}), takes
the following simplified form: for any $0<\varepsilon'< {1 \over
34}$, there exists $c_1>0$ such that
 \begin{equation}\label{rwrers:hirscheq}
    P\{ \max_{0 \le x \le N} V(x) < N^{{1 \over 2}-\varepsilon'}
    \} \sim c_{1} N^{-\varepsilon'}, \qquad N \to \infty.
 \end{equation}

\medskip

\subsection{Quenched results}
\label{rwrers:subs:excr}

 We define, for any $x\in \z$,
\begin{eqnarray*}
    \tau(x) := \min\left\{ n\ge 1: \, X_n =x \right\}, \qquad \min
    \emptyset := \infty.
    \label{rwrers:tau}
\end{eqnarray*}

\noindent (Note in particular that when $X_0=x$, then $\tau(x)$ is
the first {\it return} time to $x$.) Throughout the paper, we write
$P_\omega^x \{\cdot\} := P_\omega \{ \, \cdot \; | \, X_0=x\}$ (thus
$P_\omega^0 = P_\omega$) and denote by $E_\omega^x$ the expectation
with respect to $P_\omega^x$.

Recalling that $\omega_i/(1-\omega_i)=\ee^{-(V(i)-V(i-1))},$ we get,
for any $r<x<s$,
\begin{equation}
    P_\omega^x\{ \tau(r)< \tau(s)\} =
    \sum_{j=x}^{s-1} \ee^{V(j)} \left( \, \sum_{j=r}^{s-1} \ee^{V(j)}
    \right) ^{\! \! -1}.
    \label{rwrers:zeitouni}
\end{equation}

\noindent This result is proved in \cite{zeitouni}, see formula
(2.1.4).

 The next result, which gives a simple bound for the
expectation of $\tau(r) \wedge \tau(s)$ when the walk starts from a
site $x\in (r,s)$, is essentially contained in Golosov
\cite{golosov}; its proof can be found in \cite{shi-zindy}. For any
integers $r<s$, we have
 \begin{equation}
    \max_{x\in (r,\, s) \cap \z} E_\omega^x[ \tau(s) \, {\bf
    1}_{ \{ \tau(s) < \tau(r)\} }] \le (s-r)^2 \exp \left[
    \max_{r\le i \le j \le s} (V(j)-V(i)) \right].
    \label{rwrers:golosov2}
 \end{equation}

We will also use the following estimate borrowed from Lemma 7 of
Golosov \cite{golosov}: for $\ell \ge 1$ and $x<y$,
\begin{equation}
    P_\omega^x \{ \tau(y) < \ell\} \le \ell \, \exp\left(
    - \max_{x\le i<y} [V(y-1) -V(i)] \right) .
    \label{rwrers:golosov-valley}
\end{equation}

\noindent Looking at the environment backwards, we get: for $\ell
\ge 1$ and $w<x$,
\begin{equation}
    P_\omega^x \{ \tau(w) < \ell\} \le \ell\, \exp\left(
    - \max_{w<i\le x} [V(w+1) -V(i)] \right) .
    \label{rwrers:golosov-valley-backwards}
\end{equation}

Finally we quote an important result about excursions of Sinai's
walk (for detailed discussions, see Section 3 of \cite{dembo-gantert-peres-shi}). Let $b\in \z$
 and $x\in \z$, and consider $L(\tau(b), x)$ under $P_\omega^b$. In words,
 we look at the number of visits to the
site $x$ by the random walk (starting from $b$) until the first
return to $b$. Then there exist constants $c_{2}$ and $c_{3}$ such
that
\begin{equation}
    c_{2} \, \ee^{-[V(x)-V(b)]} \le E_\omega^b [ L(\tau(b), x)]
    \le c_{3} \, \ee^{-[V(x)-V(b)]} .
    \label{rwrers:E(excr)}
\end{equation}

\medskip

\section{Good environment-scenery and proof of Theorem \ref{rwrers:t:main}}
\label{rwrers:s:stepA}

For any $j \in \n^*$, we define
\begin{eqnarray*}
    d^+(j)
 &:=& \min \left\{ n\ge 0: \; V(n) \ge j\right\},
    \label{rwrers:d}
    \\
    b^+(j)
 &:=& \min \left\{ n\ge 0: \; V(n) = \min_{0\le x\le d^+(j)}
    V(x) \right\}.
    \label{rwrers:b}
\end{eqnarray*}

\noindent These variables enable us to consider the valley
$(0,b^+(j),d^+(j))$. Similarly, we define
\begin{eqnarray*}
    d^-(j)
 &:=& \max \left\{ n\le 0: \; V(n) \ge j\right\},
    \label{rwrers:dd}
    \\
    b^-(j)
 &:=& \max \left\{ n\le 0: \; V(n) = \min_{d^-(j) \le x\le 0}
    V(x) \right\}.
    \label{rwrers:bb}
\end{eqnarray*}

\noindent In the next sections, we will be frequently using the
following elementary estimates.
\begin{lemma}
 \label{rwrers:abc}
For any $\varepsilon'>0$, we have, $P$-almost surely for all large
$j$,
\begin{eqnarray*}
j^{2-\varepsilon'} \le |b^\pm(j)| < |d^\pm (j)| \le
    j^{2+\varepsilon'}.
\end{eqnarray*}
\end{lemma}

\medskip

\noindent {\it Proof.} Fix $\varepsilon'>0.$ Let us consider the
sequence $(j_p)_{p \ge 1}$ defined by $j_p:=p^{12/\varepsilon'}$ for
all $p \ge 1.$ Using (\ref{rwrers:hirscheq}), we obtain $
\sum_{p\ge1} P\{d^+ (j_p) > \frac{1}{3} j_p^{2+\varepsilon'}\}<
\infty.$ Therefore, Borel-Cantelli lemma implies that, $P$-almost
surely, $d^+ (j_p) \le \frac{1}{3} j_p^{2+\varepsilon'}$ for all
large $p,$ say $p \ge p_0.$ We fix a realization of $\omega$ and
consider $j_p\le j \le j_{p+1}$ with $p \ge p_0.$ Since $d^+ (j) \le
d^+ (j_{p+1}),$ we get
$$d^+ (j) \le \frac{1}{3} j_{p+1}^{2+\varepsilon'} \le
j^{2+\varepsilon'} \frac{1}{3}
\left(\frac{j_{p+1}}{j}\right)^{2+\varepsilon'} \le
j^{2+\varepsilon'} \frac{1}{3}
\left(\frac{j_{p+1}}{j_p}\right)^{2+\varepsilon'}=
j^{2+\varepsilon'} \frac{1}{3}
(1+p^{-1})^{\frac{12(2+\varepsilon')}{\varepsilon'}},
$$

\noindent which yields $d^+ (j)\le j^{2+\varepsilon'}$ for all large
$j.$ In a similar way, we can prove that $j^{2-\varepsilon'} \le
\nu^+ ((-\infty,-j^{1-\kappa}]) \le d^+ (j)$ for some $\kappa>0$ and
all large $j,$ which implies $j^{2-\varepsilon'} \le b^+ (j)$ for
all large $j.$ Moreover, the arguments are the same to prove that,
$P$-almost surely, $j^{2-\varepsilon'} \le |b^-(j)| < |d^- (j)| \le
j^{2+\varepsilon'}$ for all large $j.$
 \hfill$\Box$

\bigskip

 To introduce the announced ``good" environment-scenery, we fix $\varepsilon>0$
such that assumption of Part {\it (i)} of Theorem
\ref{rwrers:t:main} holds. For $\alpha \in (0,1)$ (which will depend
on $\varepsilon$), $0<c_{4}<1/6$, and  $j \ge 100$, we define
\begin{eqnarray}
    \gamma_{0}(j)&:=& j,
    \label{rwrers:gamma0}
    \nonumber
    \\
    \gamma_{i}(j)&:=& j^{(1-\alpha)^i} = (\gamma_{i-1}(j))^{1-\alpha} , \qquad i
 \ge 1,
    \label{rwrers:gamma}
    \\
    \varepsilon_{i}(j)&:=& \exp \bigg\{ -c_{4} \gamma_{i+2}(j) \bigg\}, \qquad i \ge 0.
    \label{rwrers:epsilon}
\end{eqnarray}

\noindent For convenience of notation we define
$\varepsilon_{-1}(j):=\varepsilon_{0}(j).$ In words,
$(\gamma_{i}(j))_{i \ge 0}$ represents a decreasing sequence of
distances, which enables us to classify the sites according to the
value of $V(x)-V(b^+(j))$.

Write $\log_p$ for the $p$-th iterative logarithmic function. Fix
$\varepsilon':= \min \{1/ 35,\varepsilon / 2\}>0,$ and introduce,
for $j \ge 100$,
\begin{equation}\label{rwrers:kappa}
    M(j):=\inf\left\{n \ge 0: \gamma_{n}(j) \le (\log_2 j)^{1-\alpha
     \over 2+ \varepsilon'} \right\}.
\end{equation}

\noindent By definition of $M(j),$ we have
 \begin{eqnarray*}
 \label{rwrers:gammakappain}
    \gamma_{M(j)-1}(j) \in \left[(\log_2 j)^{1-\alpha \over 2+ \varepsilon'},
    (\log_2 j)^{1 \over 2+ \varepsilon'}\right].
 \end{eqnarray*}

\noindent Moreover, in view of (\ref{rwrers:gamma}) and since
$\gamma_{M}(j)$ belongs to $[(\log_2 j)^{(1-\alpha)^2 \over 2+
\varepsilon'}, (\log_2 j)^{1-\alpha\over 2+ \varepsilon'}],$ we get
that
\begin{equation}
\label{rwrers:kappaequal} M(j) \sim {1 \over |\log (1-\alpha)|} \,
\log_2 j, \qquad j \to \infty.
% M=\left\lfloor {\log_2 j - \log_4j + \log
%(2+\varepsilon')-\log(1-\alpha)
%\over |\log (1-\alpha)|} \right\rfloor +1,
\end{equation}

\noindent Note that we choose $\alpha$ small enough such that
\begin{eqnarray}
\label{rwrers:varepsilonalpha}
  \beta&:=&(1-\alpha)^2 \, (2+ \varepsilon)- (2+\varepsilon')>0,
  \\
\label{rwrers:betaprim}
  \beta'&:=&{\varepsilon' \over 2}-\alpha>0.
\end{eqnarray}

 Then we introduce the set (the constant $c_{5}$ will be chosen small enough in (\ref{rwrers:choix:c_5}))
\begin{eqnarray*}
       \overline{\Theta}_{ M(j)-1 }(j)
 := \left[ b^+(j) - c_{5} \,  (\gamma_{M(j)-1}(j))^{2+\varepsilon'} , \, b^+(j) +
 c_{5} \,  (\gamma_{M(j)-1}(j))^{2+\varepsilon'}  \right],
    \label{rwrers:Okappa-}
    \end{eqnarray*}

\noindent and, for $i= M(j)-2,...,1,0,$ the sets (the constant
$c_{6} \ge 1$ will be chosen large enough in
(\ref{rwrers:choix:c_6}))
    \begin{eqnarray*}
       \overline{\Theta}_{i}(j)
 :=\left[ b^+(j) - c_{6} \,  (\gamma_{i}(j))^{2+\varepsilon'} , \, b^+(j) +
 c_{6} \,  (\gamma_{i}(j))^{2+\varepsilon'}  \right]
 \setminus \bigcup_{p=i+1}^{M(j)-1}   \overline{\Theta}_{p}(j).
    \label{rwrers:O}
\end{eqnarray*}

\noindent Observe that the sets $(\overline{\Theta}_{i}(j))_{0 \le i
\le M(j)-1}$ form a partition of the interval $[ b^+(j) - c_{6} \,
j^{2+\varepsilon'} , \, b^+(j) + c_{6} \, j^{2+\varepsilon'}].$ The
final sets we consider are given, for $0 \le i \le M(j)-1$, by
\begin{eqnarray*}\label{rwrers:Thetabar}
  \Theta_{i}(j):=  \overline{\Theta}_{i}(j) \cap I(j),
\end{eqnarray*}

\noindent where $I(j):=[\nu^+((-\infty,-j]),d^+(j)]$. Note that
$\nu^+((-\infty,-j])<d^+(j)$ on $A(j)$ which will be defined in
(\ref{rwrers:E+S}). In this case, the sets $(\Theta_{i}(j))_{0 \le i
\le M(j)-1}$ form a partition of $I(j)$ into annuli (since $c_6 \ge
1$). Loosely speaking, the set $\Theta_{i}(j)$ contains the sites
$x$ satisfying $V(x)-V(b^+(j)) \approx \gamma_i(j).$ To cover
$[d^-(j), d^+(j)]$, we define
\begin{equation}\label{rwrers:Thetabar-}
 \Theta_{-1}(j):=\left[- j^{2+\varepsilon'},
    j^{2+\varepsilon'}\right] \cap
   \left[d^-(j),  \nu^+((-\infty,-j])\right].
\end{equation}

 Moreover, for the environment on $\z^+$, we introduce the events
\begin{eqnarray}
    A^{env}_{1}(j)
 &:=& \left\{ - 4j \le V(b^+(j)) \le - 3j \right\},
    \label{rwrers:V(b)}
    \\
    A^{env}_{2}(j)
 &:=& \left\{ \max_{0\le x \le y \le b^+(j)}
    [V(y)-V(x)] \le {j \over 4} \right\}.
    \label{rwrers:V(y)-V(x)}
\end{eqnarray}

\noindent The first event ensures that the valley considered is
``deep enough" and the second one that the particle reaches the
bottom of the valley ``fast enough". To control the time spent by
the particle in different $\Theta_{i}(j)$ during an excursion from
$b^+(j)$ to $b^+(j)$, we define
\begin{eqnarray}
    A^{env}_{ann}(j)
 &:=& \bigcap_{i=0}^{M(j) -2} \left\{ \sum_{x \in \Theta_{i}(j)} \ee^{- [V(x)-V(b^+(j))]} \le
    (\varepsilon_{i}(j))^2  \right\}=:\bigcap_{i=0}^{M(j) -2}
    A^{env}_{ann,i}(j).
    \label{rwrers:LOi}
\end{eqnarray}

For the environment on $\z^-,$ let
\begin{eqnarray}
    B^{env}(j)
 &:=& \left\{ V(b^-(j)) \le -{j \over 6} \, , \, \max_{d^-(j)\le x \le y \le
    0} [V(y)-V(x)] \le {j \over 3} \right\},
    \label{rwrers:V(x)-V(y):-}
\end{eqnarray}

\noindent which ensures that the particle will not spent too much
time on $\z^-.$

 Recalling that $\xi^-_x= \max \{-\xi_x, 0 \}$, we define for the
scenery
\begin{eqnarray}
       A^{sce}_{i}(j)
 &:=& \left\{ \max_{x \in \Theta_{i}(j)}
    \xi^-_x < (\varepsilon_{i}(j))^{-1/2} \right\},  \qquad -1 \le i \le
   M(j)-2,
    \label{rwrers:Sp}
\end{eqnarray}

\noindent which ensures that the scenery does not reach excessive
negative value in each $\Theta_{i}(j)$. In order to force the
scenery in a neighborhood of the bottom (where the particle is
concentrated), to be close to $a=\mathrm{ess } \sup \xi_0$, we fix
$\rho \in (0,1)$ and introduce
\begin{eqnarray}
     A^{sce}_{M(j)-1}(j)
 &:=& \left\{ \min_{x \in \Theta_{M(j)-1}(j)}
    \xi_x \ge a - \rho \right\}.
    \label{rwrers:Skappa}
\end{eqnarray}

 We set
\begin{eqnarray*}
\nonumber
    A^{env}(j) := A^{env}_{1}(j) \cap A^{env}_{2}(j) \cap A^{env}_{ann}(j),
     \qquad A^{sce}(j) := \bigcap_{i=-1}^{M(j)-1} A^{sce}_{i}(j).
\end{eqnarray*}

\noindent Moreover, we define
\begin{equation}
\label{rwrers:E+S}
   A(j) := A^{env}(j)  \cap B^{env}(j) \cap A^{sce}(j).
\end{equation}

\noindent A pair $(\omega,\xi)$ is a ``good" environment-scenery if
$(\omega,\xi) \in  A(j)$ for infinitely many  $j\in \n$.

For future use, let us note that for $\omega \in B^{env}(j) \cap
A^{env}_{2}(j)$, we have
\begin{equation}
    \max_{d^-(j) \le x\le y\le b^+(j)} [V(y) - V(x)] \le
    {2j \over 3}.
    \label{rwrers:maxV(b-,d+)}
\end{equation}

To prove Theorem \ref{rwrers:t:main}, we need two propositions,
whose proofs are respectively postponed until Sections
\ref{rwrers:proofenv} and \ref{rwrers:proofpart}. The first one
ensures that almost all pair $(\omega,\xi)$ is a ``good''
environment-scenery, while the second one describes the behavior of
the particle in a ``good'' environment.

\medskip

\begin{proposition}
 \label{rwrers:p:good}
Under assumptions
$(\ref{rwrers:ellipticity})$--$(\ref{rwrers:Var>0})$, we have that
$P \otimes Q$-almost all $(\omega,\xi)$ is a ``good"
environment-scenery. More precisely, $P \otimes Q$-almost surely,
there exists a random sequence $(m_k)_{k \ge 1}$ such that $m_k \ge
k^{3k}$ and $(\omega,\xi)$ is a good environment-scenery along
$(m_k)_{k \ge 1}$, i.e., $(\omega,\xi) \in  A(m_k)$, for all $k\ge
1$.
\end{proposition}

\medskip

In fact $(m_k)_{k \ge 1}$ is constructed in the following way. Let
us first introduce the sequence $j_p:=p^{3p}$ for $p \ge 0.$ We
define then $(m_k)_{k \ge 1}$ by $m_1:=\inf\{j_p \ge 0: (\omega,\xi)
\in A(j_p)\}$ and $m_k:=\inf\{j_p > m_{k-1}: (\omega,\xi) \in
A(j_p)\}$ for $k \ge 2.$ Then, Proposition \ref{rwrers:p:good} means
that $m_k \to \infty,$ $k \to \infty,$ $P \otimes Q$-almost surely.
Before establishing the proposition about the behavior of the
particle, we extract a random sequence $(n_k)_{k \ge 1}$ from
$(m_k)_{k \ge 1}$ such that
\begin{eqnarray}
    \sum_{k \ge 1} \varepsilon_{M(n_k)}(n_k) < \infty .
    \label{rwrers:sum}
\end{eqnarray}

\noindent In fact, we consider the random sequence defined by
$n_1:=\inf\{m_p \ge 1: \varepsilon_{M(m_p)}(m_p) \le 1\}$ and
$n_k:=\inf\{m_p > n_{k-1}: \varepsilon_{M(m_p)}(m_p) \le
\frac{1}{k^2}\}$ for $k \ge 2.$

To ease notations, we write throughout the paper, $d_k^+ :=
d^+(n_k)$, $\tau_k^+ :=\tau(d_k^+)$, $b_k^+ := b^+(n_k)$ and $d_k^-
:= d^-(n_k)$, $\tau_k^- :=\tau(d_k^-).$ Moreover, we define, for all
$k \ge 1$,
\begin{equation}
     t_k := \lfloor \ee^{n_k} \rfloor .
    \label{rwrers:nk}
\end{equation}

\medskip

\begin{proposition}
 \label{rwrers:P}
For $P \otimes Q$ almost all $(\omega,\xi)$, we have that,
$P_{\omega}$-a.s., for all large $k,$
\begin{eqnarray}
 L \left(t_k, \Theta_{-1}(n_k) \right) &\le& \varepsilon_{-1}(n_k)
    \, t_k,
    \label{rwrers:P0}
    \\
    L \left(t_k, \Theta_{i}(n_k) \right) &\le& \varepsilon_{i}(n_k)
    \, t_k, \qquad 0 \le i \le M(n_k)-2 ,
    \label{rwrers:P1}
    \\
    \tau_k^+ \wedge \tau_k^- &>&
    t_k.
    \label{rwrers:P2}
\end{eqnarray}
\end{proposition}

\begin{remark}
There is no measurability problem for the events described in
Proposition \ref{rwrers:P}, see the beginning of Section
\ref{rwrers:proofpart}. Similar arguments apply to the forthcoming
events.
\end{remark}

\medskip

\noindent {\it Proof of Theorem \ref{rwrers:t:main}.}

\noindent {\it Proof of Part (i).} For any $\delta>0,$ we define
$\varepsilon^{(\delta)}(j) := \sum_{i=-1}^{M(j)-2}
\varepsilon_{i}^{\delta}(j).$ Recalling (\ref{rwrers:Zdecompo}), we
use Proposition \ref{rwrers:P} and Lemma \ref{rwrers:abc} to obtain,
for $\p \otimes Q$-almost all realization of $\omega$, $\xi$ and
$(X_j)_{j\ge0}$,
\begin{eqnarray}
  \sum_{j=0}^{t_k} \xi(X_j) &\ge&
  (1-\varepsilon^{(1)}(n_k))\, t_k \, \bigg( \min_{x \in \Theta_{M(n_k)-1}(n_k)}
    \xi_x \bigg)- \sum_{i=-1}^{M(n_k)-2}
  \varepsilon_i(n_k) \, t_k \, \left( \max_{x \in \Theta_i(n_k)}  \xi_x^-
  \right),
  \nonumber
\end{eqnarray}

\noindent for all large $k.$ Then, Proposition \ref{rwrers:p:good}
implies
\begin{eqnarray}
  \sum_{j=0}^{t_k} \xi(X_j)&\ge& (1-\varepsilon^{(1)}(n_k))\, t_k \, (a-\rho)- \sum_{i=-1}^
  {M(n_k)-2} \sqrt{\varepsilon_i(n_k)} \, t_k
  \nonumber
   \\
  &\ge& (1-\varepsilon^{(1)}(n_k))\, t_k \, (a-\rho)- {\varepsilon^{(1/2)}(n_k)} \, t_k,
  \label{rwrers:inequprop}
\end{eqnarray}

\noindent for all large $k$. We claim that, for any $\delta>0$ and all large $j$,
\begin{eqnarray}
   \varepsilon^{(\delta)}(j) \le \sum_{i=-1}^{ M(j)}
   \varepsilon_{i}^\delta(j) \le 2\left(1+ {1 \over \delta}\right)
   \varepsilon_{M(j)}^\delta(j).
   \label{rwrers:lemmemodif}
\end{eqnarray}

\noindent To prove (\ref{rwrers:lemmemodif}), we observe that
\begin{eqnarray*}
  \sum_{i=-1}^{M(j)} \varepsilon_{i}^\delta(j) \le  2 \, \varepsilon_{M(j)}^\delta(j)+
  \sum_{i=0}^{M(j)-1}  \int_{ \varepsilon_{i}(j)}^{ \varepsilon_{i+1}(j)} \hspace{-0.3cm}
{ \varepsilon_{i}^\delta(j) \over
    \varepsilon_{i+1}(j)- \varepsilon_{i}(j)} \d x.
\end{eqnarray*}
\noindent Recalling (\ref{rwrers:epsilon}), we have that
$\varepsilon_{i+1}(j)- \varepsilon_{i}(j)
=\varepsilon_{i+1}(j)\left(1-\ee^{-c_{4}(\gamma_{i+2}(j)-\gamma_{i+3}(j))}\right).$
Recalling  (\ref{rwrers:gamma}) we get that
$\gamma_{i+2}(j)-\gamma_{i+3}(j)=\gamma_{i+2}(j)(1-\gamma_{i+2}^{-\alpha}(j))$.
Since (\ref{rwrers:gamma}) and (\ref{rwrers:kappa}) imply
$\gamma_{i+2}(j) \ge \gamma_{M(j)+2}(j) \ge (\log_2 j)^{(1-\alpha)^4
\over 2+ \varepsilon'}$ for $0 \le i \le M(j)$, we obtain that
$\gamma_{i+2}(j)-\gamma_{i+3}(j) \ge \gamma_{i+2}(j)/2$, for all
large $j$ and for $0 \le i \le M(j)$. Therefore, we get
$\varepsilon_{i+1}(j)- \varepsilon_{i}(j)
\ge\varepsilon_{i+1}(j)/2$, implying that
\begin{eqnarray*}
   \varepsilon^{(\delta)}(j) \le 2 \, \varepsilon_{M(j)}^\delta(j)+
  2 \, \sum_{i=0}^{M(j)-1}  \int_{ \varepsilon_{i}(j)}^{ \varepsilon_{i+1}(j)}
\hspace{-0.1cm} { \varepsilon_{i}^\delta(j) \over
   \varepsilon_{i+1}(j)} \d x.
\end{eqnarray*}

\noindent Furthermore, $\sum_{i=0}^{M(j)-1}  \int_{
\varepsilon_{i}(j)}^{ \varepsilon_{i+1}(j)} {
\varepsilon_{i}^\delta(j) \over \varepsilon_{i+1}(j)} \d x \le
\sum_{i=0}^{M(j)-1}  \int_{ \varepsilon_{i}(j)}^{
\varepsilon_{i+1}(j)} x^{\delta-1} \d x= \int_{
\varepsilon_{0}(j)}^{ \varepsilon_{M(j)}(j)} x^{\delta-1} \d x,$
which is less than $\varepsilon_{M(j)}^{\delta}(j) /\delta$. This
implies (\ref{rwrers:lemmemodif}).

Combining (\ref{rwrers:inequprop}) and (\ref{rwrers:lemmemodif}) and
recalling that $\varepsilon_{M(j)}^{\delta}(j) \to 0$ when $j \to
\infty$, we get
\begin{equation}
\limsup_{n \rightarrow \infty}{1 \over n} \sum_{i=0}^{n} \xi(X_i)
\ge a - \rho, \qquad \p \otimes Q\textrm{-a.s.}
 \label{rwrers:limsup}
\end{equation}

To conclude the proof, it remains only to observe that
(\ref{rwrers:limsup}) is true for all $\rho>0$ and that the
definition of $a$ implies that $\p \otimes Q $-a.s., ${1 \over n}
\sum_{i=0}^{n} \xi(X_i) \le a$, for all $n \ge 0$. \hfill$\Box$

 \bigskip

\noindent {\it Proof of Part (ii).} Using Theorem $1.5$ of
\cite{hu-shi98a}, we have that, for any $\varepsilon''>0$,
$\p$-almost surely, $\max_{0 \le i \le n} X_i \ge (\log
n)^{2-\varepsilon''}+1$, for all large $n$. This implies
\begin{equation}
\label{rwrers:neg1}
 \sum_{i=0}^{n} \xi(X_i) \le a \, n - \max_{0 \le
x \le \lceil (\log n)^{2-\varepsilon''} \rceil} \xi^-_x.
\end{equation}

\noindent By assumption, there exists $\varepsilon>0$ such that
$Q\left\{\xi_0 <- \lambda \right\} \ge (\log
\lambda)^{-2+\varepsilon}.$ Therefore, fixing
$\varepsilon''<\varepsilon,$ we get for $k\ge1$ and all $N \ge 1,$
\begin{equation}\label{rwrers:neg2}
   Q\left\{ \max_{0 \le
x \le N} \xi^-_x < k \, a \, \ee^{N^{1 \over 2-\varepsilon''}}
\right\} \le \exp\left\{-c_{7} N^{\delta}\right\},
\end{equation}

\noindent where $\delta:=1- {2-\varepsilon \over
2-\varepsilon''}>0.$

 We choose $N_p:=\lfloor (\log p)^{T} \rfloor$ for $p\ge 1$
with $T$ large enough such that $T \delta>1$. Therefore,
(\ref{rwrers:neg2}) and the Borel--Cantelli lemma imply that,
$Q$-almost surely, there exists $p_0(\xi)$ such that
\begin{equation}\label{rwrers:neg3}
   \max_{0 \le
x \le N_p} \xi^-_x \ge k \, a \, \ee^{N_p^{1 \over
2-\varepsilon''}},
\end{equation}

\noindent for $p\ge p_0(\xi).$ Fixing a realization of $\xi$, we
define $p(n)$ by
\begin{equation}\label{rwrers:pn}
  N_{p(n)} \le  \lceil (\log n)^{2-\varepsilon''} \rceil \le
  N_{p(n)+1},
\end{equation}

\noindent for all $n$ such that $p(n)\ge p_0(\xi).$ This yields
$$
\max_{0 \le x \le \lceil (\log n)^{2-\varepsilon''} \rceil} \xi^-_x
\ge \max_{0 \le x \le N_{p(n)}} \xi^-_x \ge k \, a \,
\ee^{N_{p(n)}^{1 \over 2-\varepsilon''}},
$$

\noindent the last inequality being a consequence of
(\ref{rwrers:neg3}). Therefore, we obtain

\begin{eqnarray}
\nonumber
  \max_{0 \le x \le \lceil (\log n)^{2-\varepsilon''} \rceil} \xi^-_x &\ge& k  a
  \exp \left\{ {\lceil (\log n)^{2-\varepsilon''} \rceil}^{1
\over 2-\varepsilon''}\right\} \, \exp \left\{- \left( {\lceil (\log
n)^{2-\varepsilon''} \rceil}^{1 \over 2-\varepsilon''}-N_{p(n)}^{1
\over 2-\varepsilon''} \right) \right\}
 \\
 \nonumber
   &\ge& k  a  n  \exp \left\{- \left( N_{p(n)+1}^{1 \over 2-\varepsilon''}-N_{p(n)}^{1
\over 2-\varepsilon''} \right) \right\},
\end{eqnarray}

\noindent the second inequality being a consequence of
(\ref{rwrers:pn}). Moreover, we easily get that $N_{p(n)+1}^{1 \over
2-\varepsilon''}-N_{p(n)}^{1 \over 2-\varepsilon''}\to 0$, when $n
\to \infty$, implying that for all large $n$,
\begin{equation}\label{rwrers:neg4}
    \max_{0 \le x \le \lceil (\log n)^{2-\varepsilon''} \rceil} \xi^-_x
\ge  {k \over 2} \, a n.
\end{equation}

\noindent Assembling (\ref{rwrers:neg1}) and (\ref{rwrers:neg4}), we
get that $\p \otimes Q$-almost surely,
 $\limsup_{n \to \infty} {1 \over n} \, \sum_{i=0}^{n} \xi(X_i) \le a \, (1-{k \over 2}).$
We conclude the proof by sending $k$ to infinity.
 \hfill $\Box$
\medskip

\begin{remark}\label{rwrers:r:negativpart}
 It is possible to give more precision in the case $(ii).$  Indeed,
using the same arguments, we can prove that if $Q\{\xi^-_0 >
\lambda\} \ge {1 \over (\log \lambda)^{\alpha}},$ for some $\alpha <
2,$ then we have, for any $\varepsilon'>0,$ that $\lim_{n
\rightarrow \infty} n^{-\frac{2}{\alpha}+\varepsilon'}  Z_n
=-\infty,$ $\p \otimes Q$-almost surely.
\end{remark}

\section{Proof of Proposition \ref{rwrers:P}}
   \label{rwrers:proofpart}

Let us first explain why the events described in Proposition
\ref{rwrers:P} (more precisely in (\ref{rwrers:P0}) and
(\ref{rwrers:P1})) are measurable. Since the sequences
$(n_k)_{k\ge0}$ and $(t_k)_{k\ge0}$ can be explicitly constructed,
 $\omega \mapsto (n_k)_{k\ge0}(\omega)$ and $\omega \mapsto
(t_k)_{k\ge0}(\omega)$ are measurable. Moreover, this implies that
$\Theta_{i}(n_k)$ is measurable, for any $-1 \le i \le M(n_k)-2.$
Now, let us write
$$
L(t_k,\Theta_{i}(n_k))= \sum_{x \in \z} L(t_k,x) {\bf 1}_{\{x \in
\Theta_{i}(n_k) \}}.
$$
Since the $\Theta_{i}(n_k)$'s are measurable, so are the random
variables ${\bf 1}_{\{x \in \Theta_{i}(n_k) \}}.$ To the other hand,
the measurability of $L(t_k,x),$ for any $x \in \z,$ is obvious,
being the composition of the measurable applications  $\omega
\mapsto (t_k)_{k\ge0}(\omega)$ and $t \mapsto L(t,x).$

We now proceed to the proof of Proposition \ref{rwrers:P}.  To get
(\ref{rwrers:P2}), we observe that
$$P_{\omega} \left\{\tau_k^+ \wedge \tau_k^- \le
    t_k \right\} \le P_{\omega} \left\{\tau_k^+ \le
    t_k \right\}+P_{\omega} \left\{\tau_k^- \le
    t_k \right\}.
$$

\noindent Then using (\ref{rwrers:golosov-valley}),
(\ref{rwrers:golosov-valley-backwards}) and (\ref{rwrers:V(b)}),
(\ref{rwrers:V(x)-V(y):-}) we obtain
\begin{eqnarray*}
\label{rwrers:exittime} P_{\omega} \left\{\tau_k^+ \wedge \tau_k^-
\le
    t_k \right\}\le t_k \, (\ee^{-4n_k} + \ee^{-7n_k /6}) \le 2 \, \ee^{-n_k
    /6},
\end{eqnarray*}

\noindent Since $n_k \ge k$, this yields
$$
\sum_{k\ge0} P_{\omega} \left\{\tau_k^+ \wedge \tau_k^- \le
    t_k \right\} \le 2 \, \sum_{k\ge0} \ee^{-n_k
    /6}< \infty.
$$

 \noindent We conclude by using the Borel--Cantelli lemma.

To prove (\ref{rwrers:P1}), we apply the strong Markov property at
$\tau(b_k^+)$ and get for $0 \le i \le M(n_k)-2 $,
\begin{eqnarray*}
 &&P_{\omega} \left\{L\left(t_k, \Theta_{i}(n_k) \right) \ge
 \varepsilon_{i}(n_k) \, t_k \right\}
 \\
 &\le& P_{\omega}^{b_k^+} \left\{L\left(t_k, \Theta_{i}(n_k) \right) \ge
 \varepsilon_{i}(n_k) \, t_k - \lambda_k \right\} + P_{\omega} \left\{\lambda_k
  \le \tau(b_k^+) \le \tau_k^- \right\} + P_{\omega} \left\{\tau_k^- \le
  \tau(b_k^+)\right\},
 \end{eqnarray*}

\noindent for any $\lambda_k \ge 0.$ By (\ref{rwrers:zeitouni}),
(\ref{rwrers:V(y)-V(x)}) and Lemma \ref{rwrers:abc}, we get, for all
large $k,$
$$P_{\omega} \left\{\tau_k^- \le \tau(b_k^+)\right\} \le {b_k^+ \ee^{n_k / 3} \over
\ee^{n_k}} \le n_k^{2+\varepsilon'} \ee^{-2n_k / 3} .
$$

\noindent Since $P_{\omega} \left\{\lambda_k \le \tau(b_k^+) \le
\tau_k^- \right\} \le \lambda_k^{-1} E_\omega\left[ \tau(b_k^+) \,
{\bf 1}_{ \{\tau(b_k^+) \le \tau_k^- \} } \right]$,
(\ref{rwrers:golosov2}) and (\ref{rwrers:maxV(b-,d+)}) yield
$$
P_{\omega} \left\{\lambda_k \le \tau(b_k^+) \le \tau_k^-
    \right\} \le {(b_k^+-d_k^-)^2 \over \lambda_k} \; \ee^{2n_k / 3} \le
    {2n_k^{2 \, (2+\varepsilon')} \over \lambda_k} \; \ee^{2n_k /
    3},
$$

\noindent for all large $k,$ the second inequality being a
consequence of Lemma \ref{rwrers:abc}. Choosing
$\lambda_k:=\ee^{5n_k / 6}$, we obtain, for all large $k$,

\begin{equation}\label{rwrers:P13}
    P_{\omega} \left\{\lambda_k \le \tau(b_k^+) \le \tau_k^-
    \right\} + P_{\omega} \left\{\tau_k^- \le \tau(b_k^+)\right\} \le \ee^{- n_k / 7}.
\end{equation}

\noindent To treat $P_{k,i}:=P_{\omega}^{b_k^+} \left\{L\left(t_k,
\Theta_{i}(n_k) \right) \ge \varepsilon_{i}(n_k) \, t_k - \lambda_k
\right\}$, we observe that (\ref{rwrers:nk}) implies $\lambda_k \le
2 \, \ee^{-n_k / 6}\, t_k.$ Therefore, we obtain
\begin{eqnarray*}
  P_{k,i} &\le& P_{\omega}^{b_k^+} \left\{L\left(t_k, \Theta_{i}(n_k) \right) \ge
\left(\varepsilon_{i}(n_k)- 2 \, \ee^{-n_k / 6}\right) t_k \right\}.
\end{eqnarray*}

\noindent Then, by Chebyshev's inequality, we get
\begin{eqnarray*}
  P_{k,i} &\le& {1 \over (\varepsilon_{i}(n_k)- 2 \, \ee^{-n_k / 6}) \, t_k} \, E_\omega^{b_k^+}
   \left[L\left(t_k,\Theta_{i}(n_k) \right)\right].
\end{eqnarray*}

\noindent Furthermore, observe that Sinai's walk can not make more
than $t_k$ excursions from $b_k^+$ to $b_k^+$ before $t_k.$ Since
these excursions are i.i.d., we obtain
\begin{eqnarray}
  P_{k,i}
   &\le& { t_k \over (\varepsilon_{i}(n_k)- 2 \, \ee^{-n_k / 6}) \, t_k} \, E_\omega^{b_k^+}
   \left[L\left(\tau(b_k^+),\Theta_{i}(n_k) \right)\right].
   \label{rwrers:inequPki}
   \nonumber
\end{eqnarray}

\noindent Now we recall (\ref{rwrers:E(excr)}), which implies
$E_\omega^{b_k^+}\left[L\left(\tau(b_k^+), \Theta_{i}(n_k)
\right)\right] \le c_{3} \sum_{x \in \Theta_{i}(n_k)} \ee^{-
[V(x)-V(b^+_k)]},$ for all $0 \le i \le M(n_k)-2$. Moreover, by
(\ref{rwrers:LOi}), we get for all large $k$ and for $0 \le i \le
M(n_k)-2,$
\begin{eqnarray*}
  P_{k,i} \le {c_{3} \, \left(\varepsilon_{i}(n_k) \right)^2   \over
  \left(\varepsilon_{i}(n_k) - 2 \, \ee^{-n_k /6} \right)} \le c_{8}
  \, \varepsilon_{i}(n_k),
\end{eqnarray*}

\noindent for some $c_{8}>0.$ The second inequality is a consequence
of $\varepsilon_{i}(n_k) \ge \varepsilon_{0}(n_k)$ and the fact that
$c_{4}<1/6$ implies $\ee^{-n_k /6}=o(\varepsilon_{0}(n_k))$.

 Summing from $0$ to $M(n_k)-2$ and using (\ref{rwrers:lemmemodif}),
 we get, with $c_{9}:= 2(1+ {1 \over \delta}) \, c_{8}$,
\begin{equation}
\label{rwrers:P14}
  \sum_{i=0}^{M(n_k)-2} P_{k,i} \le c_{8}
  \sum_{i=0}^{M(n_k)-2} \varepsilon_{i}(n_k) \le
   {c_{9} \, \varepsilon_{M(n_k)}(n_k)}.
\end{equation}

\noindent Assembling (\ref{rwrers:P13}), (\ref{rwrers:P14})
 and recalling (\ref{rwrers:kappaequal}), (\ref{rwrers:sum})  we obtain
$$
\sum_{k \ge 1} \sum_{i=0}^{M(n_k)-2} P_{\omega} \left\{L\left(t_k,
\Theta_{i}(n_k) \right) \ge
 \varepsilon_{i}(n_k) \, t_k \right\} \le \sum_{k \ge 1}
\left(c_{9} \, \varepsilon_{M(n_k)}(n_k)+  M(n_k) \,\ee^{- n_k /
7}\right)< \infty.
$$

\noindent This implies (\ref{rwrers:P1}) by an application of the
Borel--Cantelli lemma.

We get (\ref{rwrers:P0}) by an argument very similar to the one used
in the proof of (\ref{rwrers:P1}), the main ingredients being the
facts that $V(x)-V(b_k^+) \ge 2 \, n_k,$ for $x \in
\Theta_{-1}(n_k)$ (which is a consequence of
(\ref{rwrers:V(x)-V(y):-}), (\ref{rwrers:V(b)}) and the definition
of $d^+(j)$), and that $\Theta_{-1}(n_k)$ contains less than $2 \,
n_k^{2+\varepsilon'}$ sites (by (\ref{rwrers:Thetabar-})). We feel
free to omit the details.

\medskip

\section{Proof of Proposition \ref{rwrers:p:good}}
   \label{rwrers:proofenv}

   We now prove that, for $P \otimes Q$-almost all $(\omega,\xi)$, there
exists a sequence $(m_k)$ such that $(\omega,\xi) \in  A(m_k)$,
$\forall k\ge 1$, where $A(m_k)$ is defined in (\ref{rwrers:E+S}).

Let $j_k := k^{3k}$ ($k\ge 1$) and $\mathcal{F}_{j_{k-1}} := \sigma
\{V(x), \xi_z, \,  d^-(j_{k-1}) \le x, z \le d^+(j_{k-1}) \}$. In
the following, we ease notations by using $\gamma_i,$
$\varepsilon_i$ and $M$ instead of $\gamma_i(j_k),$
$\varepsilon_i(j_k)$ and $M(j_k).$

 If we are able to show that
\begin{eqnarray}
    \label{rwrers:6BCL}
    \sum_k P \otimes Q \left\{ A(j_k) \, | \,
    \mathcal{F}_{j_{k-1}} \right\} = \infty, \qquad \hbox{\rm $P \otimes Q$-a.s.},
    \label{rwrers:equation}
\end{eqnarray}

\noindent then L\'evy's Borel--Cantelli lemma (\cite{durrett},
p.~$237$) will tell us that $P \otimes Q$-almost surely there are
infinitely many $k$ such that $(\omega,\xi) \in A(j_k)$.

\medskip
% Recall that $|V(x)-V(x-1)| \le M= \log {1-\delta\over \delta}$ for
% any $x\in \z$.
To bound $P \otimes Q \{ A(j_k) \, | \, \mathcal{F}_{j_{k-1}} \}$
from below, we start with the trivial inequality $A(j_k) \; \supset
\; A(j_k) \cap C(j_{k-1})$, for any set $C(j_{k-1})$. We choose
$C(j_{k-1}) := C^{env}(j_{k-1}) \cap D^{env}(j_{k-1}) \cap
C^{sce}(j_{k-1})$, where
\begin{eqnarray*}
    C^{env}(j_{k-1})
 &:=& \left\{ \inf_{0 \le y \le d^+(j_{k-1})} V(y) \ge -
j_{k-1} \log^2 j_{k-1} \right\},
    \\
    D^{env}(j_{k-1})
 &:=& \left\{ \inf_{d^-(j_{k-1}) \le y \le 0} V(y) \ge -
j_{k-1} \log^2 j_{k-1} \right\},
    \\
    C^{sce}(j_{k-1})
 &:=& \left\{ \max_{ d^-(j_{k-1})\le x \le
d^+(j_{k-1})} \xi^-_x < \left(\varepsilon_{-1}(j_k)\right)^{-1/2}
\right\}.
\end{eqnarray*}

\noindent Clearly, $C(j_{k-1})$ is
$\mathcal{F}_{j_{k-1}}$-measurable. Moreover on $C^{env}(j_{k-1})
\cap A^{env}(j_k)$, we have $d^+(j_{k-1})\le \nu^+((-\infty,-j_k])
\le b^+(j_k)$.

%Recall that $E^+(j_k) = \cap_{i=0}^4 E_i^+(j_k)$ and $S(j_k) =
%\cap_{i=0}^3 S_i(j_k)$.
Let
$$
   E^{sce}_{-1}(j_k)
 := \left\{ \max_{x \in \Theta_{-1} \setminus [d^-(j_{k-1}),d^+(j_{k-1})]}
    \xi^-_x < \left(\varepsilon_{-1}(j_k)\right)^{-1/2} \right\},
$$

\noindent and consider
$$
    E^{sce}(j_k)
 := \bigcap_{i=0}^{M-1} {A}^{sce}_i(j_k) \cap
  E^{sce}_{-1}(j_k).
$$

\noindent Since $C^{sce}(j_{k-1}) \cap E^{sce}_{-1}(j_k) \subset
 {A}^{sce}_{-1}(j_k)$, it follows that
\begin{eqnarray*}
  && P \otimes Q \left\{A(j_k) \, | \,
\mathcal{F}_{j_{k-1}} \right\}
    \\
 &\ge&  P \otimes Q \left\{ P \otimes Q
\left\{ A^{env}(j_k)  \, , \, B^{env}(j_k) \, , \, E^{sce}(j_k) \, ,
\, C(j_{k-1})  \, | \, \mathcal{F}_{j_{k-1}}, \omega \right\}  \, |
\, \mathcal{F}_{j_{k-1}} \right\}
    \\
    &\ge&  P \otimes Q \left\{ {\bf
1}_{\{ A^{env}(j_k)  \, , \ B^{env}(j_k) \, , \, C(j_{k-1}) \}} P
\otimes Q \left\{ E^{sce}(j_k) \, | \, \omega
 \right\}  \, | \, \mathcal{F}_{j_{k-1}} \right\}.
\end{eqnarray*}

\noindent  Now, we suppose for the moment that we are able to prove
that there exists $c_{10}>0$ such that, for $P$-almost all $\omega$,
 \begin{equation}
  P \otimes Q
\left\{E^{sce}(j_k) \, | \, \omega
 \right\} \ge {c_{10} \over k^{1/4}}.
 \label{rwrers:PQ1}
 \end{equation}

\noindent We get
\begin{eqnarray}
  P \otimes Q \left\{A(j_k) \, | \,
\mathcal{F}_{j_{k-1}} \right\}
  &\ge& {c_{10} \over k^{1/4}} \;  P \otimes Q \left\{ A^{env}(j_k)
  \, , \ B^{env}(j_k) \, , \, C(j_{k-1}) \, | \, \mathcal{F}_{j_{k-1}} \right\}
  \nonumber
    \\
    &\ge& {c_{10} \over k^{1/4}} \;  P_k^+ \,  P_k^-  \, {\bf
1}_{C^{sce}(j_{k-1})}, \label{rwrers:minop+p-}
\end{eqnarray}

\noindent where we use the fact that $(V(x), \, x \ge 0)$ and
$(V(x), \, x < 0)$ are independent processes and introduce
\begin{eqnarray*}
    P_{k}^{+} &:=& P  \left\{ A^{env}(j_k)  \, , \ C^{env}(j_{k-1})  \, | \,
    \mathcal{F}_{j_{k-1}}\right\},
    \label{rwrers:Pplusdef}
    \\
    P_{k}^{-} &:=& P  \left\{ B^{env}(j_k)  \, , \ D^{env}(j_{k-1})    \, | \,
    \mathcal{F}_{j_{k-1}}\right\}.
    \label{rwrers:Pminusdef}
\end{eqnarray*}

\noindent To bound $P_k^+$ from below, we introduce
\begin{eqnarray*}
    E^{env}_2(j_k)
 &:=& \left\{ \max_{0 \le x \le y \le b^+(j_k)}
    [V(y)-V(x)] \le {j_k \over 4}-j_{k-1} \log^2 j_{k-1} -
    j_{k-1} - L \right\},
\end{eqnarray*}

\noindent and consider
$$
    E^{env}(j_k) :=  A^{env}_{1}(j_k)  \cap
    A^{env}_{ann}(j_k)  \cap E^{env}_2(j_k).
$$

\noindent Observe that $C^{env}(j_{k-1}) \cap \{ \max_{d^+(j_{k-1})
\le x \le y \le b^+(j_k)}
    [V(y)-V(x)] \le {j_k \over 4}-j_{k-1} \log^2 j_{k-1} -
    j_{k-1} - L \} \subset A^{env}_2(j_k)$. Thus, since
$V(d^+(j_{k-1})) \in I_{j_{k-1}} := [j_{k-1}, \,j_{k-1}+L]$, we
have, by applying the strong Markov property at $d^+(j_{k-1})$,
\begin{equation}
    P_k^+ \ge \left(\, \inf_{z \in
I_{j_{k-1}}} P_z  \left\{  E^{env}(j_k) \right\} \right) \, {\bf
1}_{C^{env}(j_{k-1})}. \label{rwrers:F+0}
\end{equation}

\noindent To bound $P_k^-$ from below, we observe the following
inclusion
$$B^{env}(j_k) \supset \left\{
\max_{d^-(j_k)\le x \le y \le d^-(j_{k-1})} [V(y)-V(x)] \le {j_k
\over 3} \right\} \cap D^{env}(j_{k-1}).$$
 \noindent Then since $V(d^-(j_{k-1}))$ belongs to $ I_{j_{k-1}}$,
the strong Markov property applied at $d^-(j_{k-1})$ yields
\begin{equation}
    P_k^- \ge \left(\, \inf_{z \in
I_{j_{k-1}}} P_z  \left\{B^{env}(j_k) \right\} \right) \, {\bf
1}_{D^{env}(j_{k-1})}. \label{rwrers:F-0}
\end{equation}

Observe that an easy calculation yields ${\bf 1}_{ C(j_{k-1})} = 1$,
$P \otimes Q$-almost surely for all large $k.$ Therefore, recalling
(\ref{rwrers:minop+p-}), (\ref{rwrers:F+0}) and (\ref{rwrers:F-0}),
the proof of (\ref{rwrers:6BCL}) boils down to showing that
\begin{eqnarray}
   \liminf_{k \to \infty} \inf_{z \in I_{j_{k-1}}} P_z \left\{ E^{env}(j_k) \right\}
    &>& 0,
    \label{rwrers:F+}
    \\
  \liminf_{k \to \infty}  \inf_{z \in I_{j_{k-1}}} P_z \left\{ B^{env}(j_k)  \right\}
    &>& 0.
    \label{rwrers:F-}
\end{eqnarray}

The rest of the section is devoted to the proof of
(\ref{rwrers:PQ1}) and (\ref{rwrers:F+}), whereas (\ref{rwrers:F-})
is an immediate consequence of Donsker's theorem.

\medskip

\subsection{Proof of (\ref{rwrers:PQ1})}
\label{rwrers:S}

\noindent Since the sets $\left\{ \Theta_{i} \right\}_{-1 \le i \le
M-1}$ are disjoint, the events $E^{sce}_{-1}(j_k)$ and
 $\left\{ A^{sce}_i(j_k) \right\}_{0 \le
i \le M -1}$ are mutually independent. We write
\begin{eqnarray*}
    P \otimes Q \left\{
E^{sce}(j_k) \, | \, \omega
 \right\}= \prod_{i=0}^{M-1}  P \otimes Q \left\{
A^{sce}_i(j_k) \, | \, \omega  \right\} \times P \otimes Q \left\{
E^{sce}_{-1}(j_k) \, | \, \omega \right\}.
 \label{rwrers:Sindepen}
\end{eqnarray*}

\noindent Thus, (\ref{rwrers:PQ1}) will be a consequence of the two
following lemmas.

\medskip

\begin{lemma}
 \label{rwrers:l:S3}
For $P$-almost all $\omega$, we have
\begin{eqnarray*}
    P \otimes Q \left\{ A^{sce}_{M-1}(j_k) \, | \, \omega
 \right\} \ge {1 \over k^{1/4}}.
    \label{rwrers:S3+}
\end{eqnarray*}
\end{lemma}

\medskip

\begin{lemma}
 \label{rwrers:l:Si}
There exists $c_{11}>0$ such that, for $P$-almost all $\omega,$
\begin{equation}
\liminf_{k \to \infty} \, \prod_{i=0}^{ M-2}  P \otimes Q \left\{
A^{sce}_{i}(j_k) \, | \, \omega \right\} \times P \otimes Q \left\{
E^{sce}_{-1}(j_k) \, | \, \omega
 \right\} \ge c_{11}.
 \label{rwrers:Si+}
\end{equation}
\end{lemma}

\bigskip

\noindent {\it Proof of Lemma \ref{rwrers:l:S3}.} Recalling
(\ref{rwrers:Skappa}), (\ref{rwrers:Okappa-}) and
(\ref{rwrers:gammakappain}), we get, $P$-almost surely,
\begin{eqnarray}
P \otimes Q \left\{ A^{sce}_{M-1}(j_k) \, | \, \omega
   \right\}\ge \exp \left\{ c_{5} \, \log q  \, \log_{2}j_k
   \right\},
   \nonumber
\end{eqnarray}

\noindent where $q:=Q \left\{\xi_0 \ge a - \rho \right\}$. Note that
the definition of $a$ implies $-\infty< \log q < 0.$ Therefore, it
remains only to observe that $\log_{2}j_k=\log k+\log_2k + \log3$
and to choose $c_{5}$ small enough such that
\begin{eqnarray} \label{rwrers:choix:c_5}
c_{5}\, \log q>-1/5,
\end{eqnarray}
\noindent to conclude the proof. \hfill$\Box$

\bigskip

\noindent {\it Proof of Lemma \ref{rwrers:l:Si}.} Recalling
(\ref{rwrers:Sp}) and that $(\xi^-_x)_{x \in \z}$ is a family of
i.i.d.\ random variables, we get, $P$-almost surely, for $0 \le i
\le M-2$,
\begin{eqnarray*}
  P \otimes Q \left\{ A^{sce}_i(j_k) \, | \, \omega
 \right\} &\ge& \Big(Q \big\{\xi^-_0 \le \varepsilon_i^{-1/2}
 \big\} \Big)^{2c_{6} \, \gamma_i^{2+\varepsilon'}}
  \nonumber
  \\
  &\ge& \exp \bigg\{2c_{6} \, \gamma_i^{2+\varepsilon'} \, \log
  \Big( 1- Q \big\{\xi^-_0 \ge \varepsilon_i^{-1/2} \big\} \Big)
   \bigg\}.
  \nonumber
\end{eqnarray*}

\noindent Then, since $Q\{\xi^-_0 \ge \varepsilon_i^{-1/2}\}$ tends
to $0$ when $k$ tends to $\infty$ and using the fact that $\log(1-x)
\ge - c_{12} \, x$ for $x \in [0,1/2)$ with $c_{12}:= 2 \log 2 >0$,
it follows that
$$
  P \otimes Q \left\{ A^{sce}_i(j_k) \, | \, \omega
 \right\} \ge \exp \bigg\{-c_{13} \,  \gamma_i^{2+\varepsilon'} \,
   Q \big\{\xi^-_0 \ge \varepsilon_i^{-1/2} \big\} \bigg\},
$$

\noindent for all large $k$, with $c_{13}:=2 \, c_{6} \, c_{12}$.
Recalling that $Q\left\{\xi^-_0 \ge \lambda \right\} \le (\log
\lambda)^{-(2+\varepsilon)}$ for $\lambda \ge \lambda_0>0$ and
(\ref{rwrers:epsilon}) we get for $k$ large enough and uniformly in
$0 \le i \le M -2$,
\begin{equation}\label{rwrers:Si++}
      P \otimes Q \left\{ A^{sce}_i(j_k) \, | \, \omega
 \right\} \ge \exp \bigg\{- c_{14} \,\gamma_i^{-\beta} \bigg\},
\end{equation}

\noindent where $\beta:=(1-\alpha)^2 \, (2+\varepsilon)-(2+\ \varepsilon')>0$
 by (\ref{rwrers:varepsilonalpha}), and
 $c_{14}:=c_{13} \left({2 \over c_{4}}\right)^{2+\varepsilon}$.
Similarly, since $E^{sce}_{-1}(j_k) \subset A^{sce}_{-1}(j_k)$ and
recalling (\ref{rwrers:Sp}), we obtain
\begin{equation}\label{rwrers:S--1}
    P \otimes Q \left\{ E^{sce}_{-1}(j_k) \, |
\, \omega \right\} \ge \exp \bigg\{- c_{15} \,\gamma_{0}^{-\beta}
\bigg\},
\end{equation}

\noindent for some $c_{15}>0.$ Combining (\ref{rwrers:Si++}) and
(\ref{rwrers:S--1}), we get
\begin{eqnarray*}
\prod_{i=0}^{M-2}  P \otimes Q \left\{ A^{sce}_i(j_k) \, | \, \omega
\right\} \times P \otimes Q \left\{ E^{sce}_{-1}(j_k) \, | \, \omega
 \right\} \ge \exp\{- c_{16} \, \sigma_{\beta}\},
\end{eqnarray*}

\noindent with $c_{16}:=\max\{c_{14},c_{15}\}$ and
$\sigma_{\beta}:=\gamma_0^{-\beta}+ \sum_{i=0}^{ M-2}
\gamma_i^{-\beta}$. By the same way we proved
(\ref{rwrers:lemmemodif}), we obtain that, for any $\beta>0,$ there
exists $c_{17} \le 1+2/\beta$ such that $\sigma_{\beta} \le c_{17}
\gamma_{M-1}^{-\beta}.$ Recalling (\ref{rwrers:gammakappain}), it
follows that $\sigma_{\beta} \to 0$ when $k \to \infty$, which
implies (\ref{rwrers:Si+}).
 \hfill$\Box$

 \medskip

\subsection{Proof of (\ref{rwrers:F+})}
\label{rwrers:F++}

To prove (\ref{rwrers:F+}), we need the following preliminary
result.
\medskip

\begin{lemma}
 \label{rwrers:sumgammap}
 For any $\delta>0$, $k\ge1$ and any $0 \le p \le M$, we have
\begin{eqnarray*}\label{rwrers:eqsumgamma}
    \sum_{i=p}^{M} \gamma_i^{\delta} \le \left( 1+ {2 \over \delta} \right)\,
    \gamma_p^{\delta}.
\end{eqnarray*}
\end{lemma}

\medskip

\noindent {\it Proof.} Observe that we easily get
\begin{eqnarray*}
 \sum_{i=p}^M \gamma_{i}^\delta \le \gamma_{p}^\delta+
  \sum_{i=p+1}^{M}  \int_{  \gamma_{i}}^{ \gamma_{i-1}} { \gamma_{i}^\delta \over
   \gamma_{i-1}- \gamma_{i} } dx.
\end{eqnarray*}

\noindent Recalling  that $\gamma_{i-1}-\gamma_{i} \ge
\gamma_{i-1}/2$, for all large $j$ and for $1 \le i \le M$, we get
\begin{eqnarray*}
 \sum_{i=p}^M \gamma_{i}^\delta \le \gamma_{p}^\delta+
  2 \, \sum_{i=p+1}^{M}  \int_{  \gamma_{i}}^{ \gamma_{i-1}} { \gamma_{i}^\delta \over
   \gamma_{i-1}} dx.
\end{eqnarray*}

\noindent Then, $ \sum_{i=p+1}^{M}  \int_{  \gamma_{i}}^{
\gamma_{i-1}} { \gamma_{i}^\delta \over  \gamma_{i-1}} dx  \le
\sum_{i=p+1}^{M} \int_{  \gamma_{i}}^{ \gamma_{i-1}} x^{\delta-1} dx
= \int_{\gamma_{M}}^{\gamma_{p}}  x^{\delta-1} dx \le
\gamma_{p}^{\delta} /\delta$ yields Lemma \ref{rwrers:sumgammap}.
\hfill$\Box$

\bigskip

We now proceed to prove (\ref{rwrers:F+}). Let
$$
a_{\ell}:=- 3j_k - \ell \, \gamma_{M}, \qquad F^{env}_1(j_k
,{\ell}):= \left\{ a_{\ell+1} \le V(b^+(j_k )) < a_{\ell} \right\}.
$$

\noindent  Denoting $\theta_k:=\lfloor \, j_k / \gamma_{M} \rfloor
-1$, the inclusion $\bigsqcup_{\ell=0}^{\theta_k} F^{env}_1(j_k
,{\ell}) \subset A_1^{env}(j_k)$ yields
\begin{eqnarray}
    \label{rwrers:PF+1}
P_z \left\{ E^{env}(j_k) \right\} &\ge& \sum_{{\ell}=0}^{\theta_k}
P_z
    \left\{F^{env}_1(j_k,{\ell}) \, , \,  A^{env}_{ann}(j_k)  \, , \, E^{env}_2(j_k)
    \right\} =: \sum_{{\ell}=0}^{\theta_k}
    P^+_{k,{\ell}}.
\end{eqnarray}

\noindent  To bound $ P^+_{k,{\ell}}$ by below for $0 \le \ell \le \theta_k$,
  we define the following levels,
\begin{eqnarray}
  \eta_{i} = \eta_{i}(j_k,\ell) &:=& a_{\ell} + \gamma_{i}, \qquad 0 \le i \le
  M,
  \label{rwrers:etai}
  \\
  \eta_{M+1} = \eta_{M+1}(j_k,\ell) &:=& a_{\ell},
  \label{rwrers:etakappa}
  \\
  \eta_{M+2} = \eta_{M+2}(j_k,\ell) &:=& a_{\ell+1}.
  \label{rwrers:etakappa+1}
\end{eqnarray}

\noindent In the following, we introduce stopping times for the
potential, which will enable us to consider a valley having ``good"
properties. Let us write
\begin{eqnarray}
  T = T(j_k,\ell) &:=& \nu^+((-\infty, \eta_{M+1}]),
   \nonumber
   \\
  \tilde{T} = \tilde{T}(j_k,\ell) &:=& \nu^+((-\infty, \eta_{M+2}]).
   \nonumber
\end{eqnarray}

\noindent Then, let us define the following stopping times, for $0
\le i \le M $,
\begin{eqnarray}
  T_i = T_i(j_k,\ell) &:=& \nu^+((-\infty, \eta_{i}]),
  \nonumber
   \\
  T'_i = T'_i(j_k,\ell) &:=& \min\left\{ n\ge T:
  \, V(n) \ge \eta_{M-i}\right\},
  \nonumber
  \\
  R_i = R_i(j_k,\ell) &:=& \min\left\{ n
     \ge T'_{i}: \, V(n) \le \eta_{M-i+1}\right\}.
  \nonumber
\end{eqnarray}

\noindent We introduce the events
\begin{eqnarray}
 G_i(j_k) &:=& \left\{ T_{i+1} - T_i \le  \gamma_{i}^{2+ \varepsilon'}
   \, , \,  \max_{T_i \le x \le y \le T_{i+1}}  [V(y)-V(x)] \le {j_k \over 5} \right\},
   \qquad 0 \le i \le M -1,
   \nonumber
   \\
 G_{M}(j_k) &:=& \left\{ T-
  T_M \le \gamma_{M}^{2+ \varepsilon'} \, , \,
   \max_{T_{M} \le x \le y \le T} [V(y)-V(x)] \le {j_k \over 5}  \right\},
   \nonumber
    \end{eqnarray}

\noindent and\begin{eqnarray}
   G'_{0}(j_k) &:=& \left\{ T'_0-
  T \le  \gamma_{M}^{2+ \varepsilon'} \, , \,  T'_0< \tilde{T}
   \right\},
   \nonumber
      \\
  G'_{i}(j_k) &:=& \left\{ T'_{i} - T'_{i-1} \le  \gamma_{M-i}^{2+ \varepsilon'}
   \, , \,  T'_{i}  < R_{i-1} \right\},  \qquad 1 \le i \le M.
   \nonumber
   \end{eqnarray}

\noindent Moreover, we set
\begin{eqnarray*}
G(j_k, \ell) := \bigcap_{i=0}^{M} G_i(j_k), \qquad  G'(j_k, \ell) :=
\bigcap_{i=0}^{M} G'_i(j_k),
\end{eqnarray*}

\noindent and
\begin{eqnarray}
  H(j_k,\ell) &:=& \left\{ \max_{0 \le x \le y \le T_0}
    [V(y)-V(x)] \le {j_k \over 5} \right\},
    \nonumber
  \\
 H'(j_k, \ell)&:=& \left\{ d^+(j_k)<R_{M}\right\}.
     \nonumber
\end{eqnarray}

\begin{figure}[htb]
\centering \setlength\unitlength{1pt}
\begin{picture}(600,300)
\put(60,15){\includegraphics[width=0.775\textwidth]{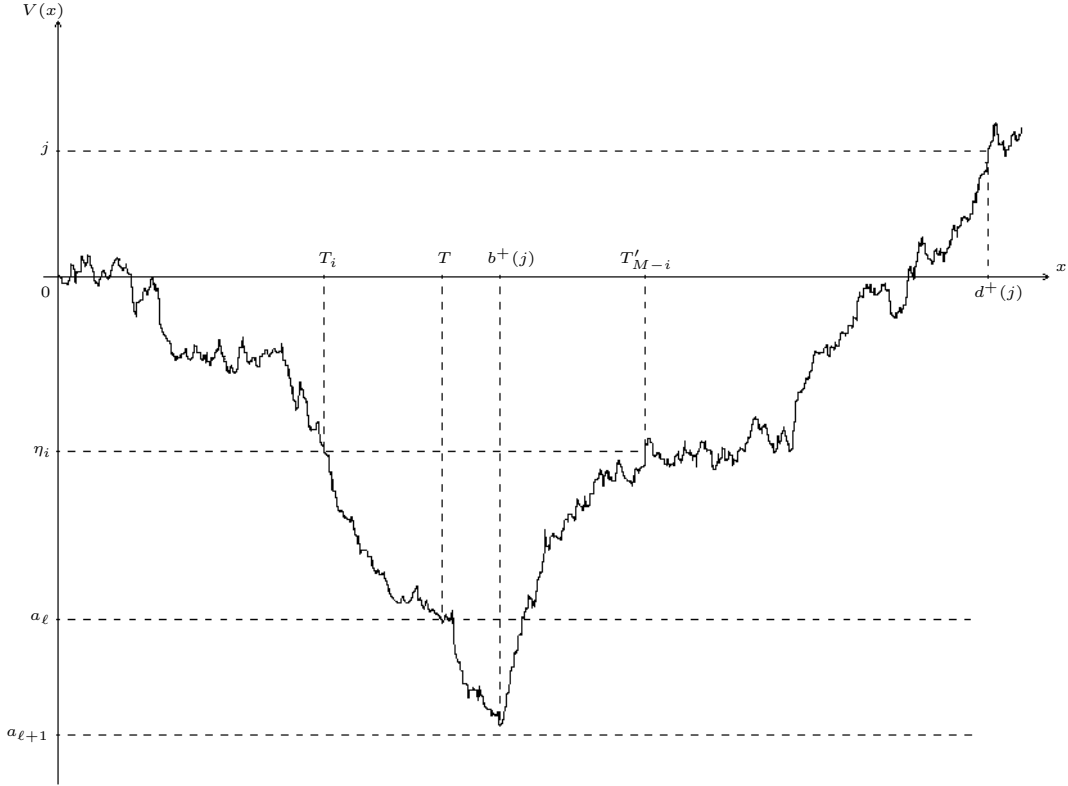}}
 \put(404,180){\tiny $d^+(j)$}
 \put(156,193){\tiny $T_i$}
 \put(270,193){\tiny $T'_{M-i}$}
\put(201,193){\tiny $T$}
 \put(220,193){\tiny $b^+(j)$}
 \put(435,190){\tiny $x$} \put(44,287){\tiny $V(x)$}
 \put(51,235){\tiny $j$}
\put(51,180){\tiny $0$}
 \put(48,121){\tiny $\eta_i$}
 \put(47,58){\tiny $a_{\ell}$}
 \put(38,14){\tiny $a_{\ell+1}$}
 \end{picture}
 \caption{Example of $\omega \in G(j_k, \ell) \cap G'(j_k, \ell) \cap H(j_k,\ell) \cap H'(j_k,\ell)$}
\end{figure}

\noindent See Figure $1$ for an example of $\omega \in G(j_k, \ell)
\cap G'(j_k, \ell) \cap H(j_k,\ell) \cap H'(j_k,\ell).$

Observe that on $G(j_k, \ell) \cap G'(j_k, \ell) \cap H'(j_k,
\ell),$ we have, for $0 \le i \le M-1,$
\begin{equation}
    [T_{i},T'_{M-i}] \supset \left\{ x \in
     [0,d^+(j_{k})]: \, V(x)-V(b^+(j_k)) \le \gamma_{i+1} \right\}.
    \label{rwrers:Tsupset}
\end{equation}

\noindent Moreover, on $G(j_k, \ell) \cap G'(j_k, \ell)$,
\begin{eqnarray*}
       T'_{M-i}-T_{i} \le 2 \sum_{p=i}^{M} \gamma_{p}
    ^{2+ \varepsilon'}, \qquad 0 \le i \le M .
   \label{rwrers:TT'}
\end{eqnarray*}

\noindent If we choose $c_{6}$ such that
\begin{eqnarray}
   \label{rwrers:choix:c_6}
c_{6} \ge 2 (1 + {2 \over 2+\varepsilon'}),
\end{eqnarray}

\noindent then Lemma \ref{rwrers:sumgammap} yields
\begin{eqnarray}
   \label{rwrers:Ri-T}
[T_{i},T'_{M-i}]
  \subset [ b^+(j_k) - c_{6} \, \gamma_{i}^{2+\varepsilon'} , \,
  b^+(j_k) + c_{6} \,  \gamma_{i}^{2+\varepsilon'}], \qquad 0 \le  i \le  M-2.
\end{eqnarray}

\noindent Recall now definition of $\Theta_i(j_k)$, so that, by
assembling (\ref{rwrers:Tsupset}) and (\ref{rwrers:Ri-T}), we have
on $G(j_k, \ell) \cap G'(j_k, \ell)\cap H'(j_k, \ell)$,
$$
 \Theta_i(j_k) \subset \left\{ x \in \z: \, V(x)-V(b^+(j_k))
 \ge \, \gamma_{i+2} \right\}, \qquad 0\le i \le M-2.
$$

\noindent An easy calculation yields $\sum_{x \in \Theta_i(j_k)}
\exp\{- [V(x)-V(b^+(j_k))]\} \le 2 c_{6} \,
\gamma_i^{2+\varepsilon'} \ee^{-\gamma_{i+2}}$, for all $0 \le  i
\le  M -2$, on  $G(j_k, \ell) \cap G'(j_k, \ell)\cap H'(j_k, \ell).$
On the other hand, since $6 \, c_{4}<1$, we get $2 c_{6} \,
\gamma_i^{2+\varepsilon'} \ee^{-\gamma_{i+2}} \le \varepsilon_i^2$,
for all large $k$ and uniformly in $0\le  i \le  M -2$. This implies
that $G(j_k, \ell) \cap G'(j_k, \ell)\cap H'(j_k, \ell) \subset
F^{env}_1(j_k,{\ell}) \cap
 A^{env}_{ann}(j_k).$ We easily observe that $G(j_k, \ell)  \cap H(j_k,\ell)
\subset  E_2^{env}(j_k)$, for all large $k$. Thus we obtain
$$
G(j_k, \ell) \cap G'(j_k, \ell) \cap H(j_k,\ell) \cap H'(j_k,\ell)
\subset F^{env}_1(j_k,{\ell}) \cap A^{env}_{ann}(j_k)  \cap
E_2^{env}(j_k).
$$

\noindent Recalling (\ref{rwrers:PF+1}), we get
\begin{eqnarray*}\label{rwrers:mino}
    P^+_{k,{\ell}} \ge P_z \left\{G(j_k, \ell) \, , \,
G'(j_k, \ell) \, , \, H(j_k,\ell) \, , \, H'(j_k,\ell)   \right\}.
\end{eqnarray*}

To bound $P_z \{G(j_k, \ell) \, , \, G'(j_k, \ell) \, , \,
H(j_k,\ell) \, , \, H'(j_k,\ell) \}$ by below, we will apply the
strong Markov property several times.

 Since $V(T'_{M}) \in I_{\eta_0}:=[\eta_0, \, \eta_0+L]$, the strong Markov
property applied at $T'_{M}$ implies, for $z \in I_{j_{k-1}}$,
\begin{eqnarray*}
    P^+_{k,\ell} \ge P_z \bigg\{ G(j_k, \ell) \, , \, G'(j_k, \ell) \, , \, H(j_k,\ell) \bigg\}
    \, \inf_{y \in I_{\eta_0}}P_y\left\{d^+(j_k) \le \nu^+((-\infty,\eta_1]) \right\} .
\end{eqnarray*}

\noindent To bound by below $P_y \{ \cdots \}$ on the right hand
side, observe that $P_y \{ \cdots \}$ is greater than  $P_{\eta_0}
\{ \cdots \}$. Moreover since $\eta_1 \ge -4j_k$ implies $j_k-
\eta_1 \le 5j_k,$ Lemma \ref{rwrers:martingale1} yields
 $$P_{\eta_0}
\{ d^+(j_k)\le \nu^+((-\infty,\eta_1])\} \ge {\eta_0 - \eta_1 \over
5j_k +L}={\eta_0 (1- \eta_0^{-\alpha}) \over 5j_k +L} \ge c_{18},$$
\noindent for all large $k$ and some $c_{18}>0$, which implies
$$
P^+_{k,\ell} \ge c_{18}  \,  P_z \left\{ G(j_k, \ell) \, , \,
G'(j_k, \ell) \, , \, H(j_k,\ell) \right\}.
$$

\noindent We now apply the strong Markov property successively at
$(T'_{M-i})_{1 \le i \le M}$ and $T$, such that
\begin{equation}
    \label{rwrers:mino2}
    P^+_{k,\ell} \ge c_{18}  \, P_z \left\{G(j_k, \ell) \, , \, H(j_k,\ell) \right\} \,
    \inf_{y \in I_{\eta_{M+1}-L}} Q'_{0,y} \, \prod_{p=1}^{M}  \inf_{y \in I_{\eta_{M-p+1}}} Q'_{p,y},
\end{equation}

\noindent where
\begin{eqnarray*}
    Q'_{p,y}&:=&P_y\left\{d^+(\eta_{M-p}) \le  \gamma_{M-p}^{2+ \varepsilon'} \, , \,
     d^+(\eta_{M-p})< \nu^+((-\infty,\eta_{M-p+2}])\right\}, \qquad
     0 \le p \le M.
\end{eqnarray*}

\noindent First, observe that $ \inf_{y \in I_{\eta_{M-p+1}}}
Q'_{p,y} \ge Q'_{p,\eta_{M-p+1}}=:Q'_p,$ for all $1 \le p \le M$ and
similarly $ \inf_{y \in I_{\eta_{M+1}-L}} Q'_{0,y} \ge
Q'_{0,\eta_{M+1}-L}:=Q'_0.$ Therefore we only have to bound from
below $Q'_p$ for $1 \le p \le M$ and $Q'_0.$ Recalling that $P\{A \,
, \, B\} \ge P\{A\}-P\{B^c\}$, we get, for $1 \le p \le M,$
\begin{eqnarray*}
\label{rwrers:mino1}
  Q'_p \ge   P_{\eta_{M-p+1}}\left\{ d^+(\eta_{M-p})<
\nu^+((-\infty,\eta_{M-p+2}])\right\}-
P_{\eta_{M-p+1}}\left\{d^+(\eta_{M-p}) \ge \gamma_{M-p}^{2+
\varepsilon'} \right\},
\end{eqnarray*}

\noindent and
\begin{eqnarray*}
\label{rwrers:mino1'}
  Q'_0 \ge   P_{\eta_{M+1}-L}\left\{ d^+(\eta_{M})<
\nu^+((-\infty,\eta_{M+2}])\right\}-
P_{\eta_{M+1}-L}\left\{d^+(\eta_{M}) \ge  \gamma_{M}^{2+
\varepsilon'} \right\}.
\end{eqnarray*}

\noindent By Lemma \ref{rwrers:martingale1}, we obtain, for $1 \le p
\le M,$
$$
P_{\eta_{M-p+1}}\left\{ d^+(\eta_{M-p})<
\nu^+((-\infty,\eta_{M-p+2}])\right\} \ge {\eta_{M-p+1}-\eta_{M-p+2}
\over \eta_{M-p}-\eta_{M-p+2} +L},
$$

\noindent and
$$
P_{\eta_{M+1}-L}\left\{ d^+(\eta_{M})<
\nu^+((-\infty,\eta_{M+2}])\right\} \ge { \eta_{M+1}-L-\eta_{M+2}
\over \eta_{M}-\eta_{M+2} +L}.
$$

\noindent Recalling (\ref{rwrers:etai}) and (\ref{rwrers:etakappa}),
we bound by below  $P_{\eta_{M-p+1}}\left\{ d^+(\eta_{M-p})<
\nu^+((-\infty,\eta_{M-p+2}])\right\}$ (for all $1 \le p \le M$) by
\begin{eqnarray*}
{\gamma_{M-p+1} \over \gamma_{M-p}}
{1-\gamma_{M-p}^{-\alpha(1-\alpha)} \over 1+L\gamma_{M-p}^{-1}}
&\ge& {\gamma_{M-p+1} \over \gamma_{M-p}}
(1-\gamma_{M-p}^{-\alpha(1-\alpha)})  (1-L\gamma_{M-p}^{-1})
\\
 &\ge& {\gamma_{M-p+1} \over \gamma_{M-p}} (1-2 \,
\gamma_{M-p}^{-\alpha(1-\alpha)}),
 \end{eqnarray*}

\noindent for all large $k.$ The first inequality is a consequence
of $(1+x)^{-1} \ge 1-x$ for any $x \in (0,1)$ and the second one is
a consequence of $0<\alpha<1.$ Similarly, recalling
(\ref{rwrers:etai}), (\ref{rwrers:etakappa}) and
(\ref{rwrers:etakappa+1}), we get, for all large $k,$
$$
P_{\eta_{M+1}-L}\left\{ d^+(\eta_{M})<
\nu^+((-\infty,\eta_{M+2}])\right\} \ge { \gamma_{M}-L \over 2 \,
\gamma_{M} +L} \ge c_{18},
$$

\noindent with $c_{18}>0.$ Moreover, combining
(\ref{rwrers:hirscheq}) and the fact that $\gamma_{M-p}\le
\gamma_{M-p}^{(2+\varepsilon')(\frac{1}{2}-\frac{\varepsilon'}{6})}$
for $0 \le p \le M$ yields
\begin{eqnarray*}
P_{\eta_{M-p+1}}\left\{d^+(\eta_{M-p}) \ge \gamma_{M-p}^{2+
\varepsilon'} \right\} &\le& c_{19} \gamma_{M-p}^{-\varepsilon'/6},
\qquad  1 \le p \le M,
\\
P_{\eta_{M+1}-L}\left\{d^+(\eta_{M}) \ge \gamma_{M}^{2+
\varepsilon'} \right\} &\le& c_{19} \gamma_{M}^{-\varepsilon'/6},
\end{eqnarray*}

\noindent for all large $k$ and for some $c_{19}>0.$ Therefore, we
obtain $Q'_0 \ge c_{20}$ for some $c_{20}>0$ and recalling
(\ref{rwrers:betaprim}) we get, for $1 \le p \le M,$
\begin{eqnarray*}
\label{rwrers:minoQp}
 Q'_p \ge {\gamma_{M-p+1} \over \gamma_{M-p}}
(1- c_{21} \gamma_{M-p}^{-\beta''}),
\end{eqnarray*}

\noindent where $\beta'':=\min\{\alpha(1-\alpha),\beta'\}>0$
($\beta'$ is defined in (\ref{rwrers:betaprim})) and $c_{21}>0.$
Observe that $\gamma_{M-p}^{-\beta''}\le \gamma_{M}^{-\beta''}$ for
$1 \le p \le M$, and that $\gamma_{M}^{-\beta''} \to 0,$ $k \to
\infty.$ Recalling the fact that $\log(1-x) \ge - c_{12} \, x,$ for
$x \in [0,1/2)$, we obtain
\begin{eqnarray*}
 \inf_{y \in I_{\eta_{M+1}-L}} Q'_{0,y} \, \prod_{p=1}^{M}  \inf_{y \in
 I_{\eta_{M-p+1}}} Q'_{p,y} \ge c_{20} {\gamma_{M} \over \gamma_0} \exp \left\{-c_{22} \,
 \sum_{p=1}^{M}
  \gamma_{M-p}^{-\beta''} \right\},
\end{eqnarray*}

\noindent where $c_{22}:=c_{12}c_{21}.$

Recall that for any $\beta''>0,$ there exists $c_{23}>0$ such that
$\sum_{p=1}^{M} \gamma_{M-p}^{-\beta''}\le c_{23} \,
\gamma_{M}^{-\beta''}.$ Then, recalling (\ref{rwrers:mino2}), this
yields, for all large $k$,
\begin{equation}
 P^+_{k,\ell} \ge c_{24} \, {\gamma_{M} \over \gamma_0}   \, P_z \left\{G(j_k, \ell)
  \, , \, H(j_k,\ell) \right\},
  \label{rwrers:mino3}
\end{equation}

\noindent with $c_{24}>0.$ To bound $P_z \{G(j_k, \ell)
  \, , \, H(j_k,\ell) \}$ from below, we apply successively the strong
  Markov property at $(T_{M-i})_{0 \le i \le M}$ such that
\begin{eqnarray*}
    P_z \left\{G(j_k, \ell) \, , \, H(j_k,\ell) \right\} \ge P_z \left\{
    H(j_k,\ell) \right\} \,
     \prod_{p=0}^{M} Q_{p},
\end{eqnarray*}

\noindent where
\begin{eqnarray*}
    Q_{p}&:=&P_{\eta_{p}}\left\{\nu^+((-\infty,\eta_{p+1}]) \le  \min\left\{d^+(\eta_{p+1}+j_k/5),
     \gamma_{p}^{2+ \varepsilon'}\right\} \right\}, \qquad 0 \le p \le M.
\end{eqnarray*}

\noindent Recall that $P\{A \, , \, B\} \ge P\{A\}-P\{B^c\}$. Then
(\ref{rwrers:hirscheq}) yields, for $1 \le p \le M,$
$$P_{\eta_p} \left\{\nu^+((-\infty,\eta_{p+1}]) \le
 \gamma_{p}^{2+ \varepsilon'} \right\} \ge1- c_{25}
\gamma_{p}^{-\varepsilon'/6},
$$

\noindent with $c_{25}>0.$ Moreover, using Lemma
\ref{rwrers:martingale1}, we get, for $1 \le p \le M,$
$$
P_{\eta_{p}}\left\{d^+(\eta_{p+1}+j_k/5) \le
 \nu^+((-\infty,\eta_{p+1}])\right\} \le c_{26} \,
{\gamma_{p} \over j_k},
$$

\noindent with $c_{26}>0.$ Therefore, observing that, for $1 \le p
\le M,$ we have ${\gamma_{p} \over j_k} \le {\gamma_{1} \over
j_k}=j_k^{-\alpha} \to 0,$ $k \to \infty,$ and using the fact that
$\log(1-x) \ge - c_{12} \, x,$ for $x \in [0,1/2)$, we get that
\begin{equation}
    \label{rwrers:mino4}
     \prod_{p=1}^{M} Q_{p} \ge \exp \left\{-c_{27} \, \sum_{p=1}^{M}
  \left(\gamma_{p}^{-\varepsilon'/6}+{\gamma_{p} \over j_k}\right)
  \right\},
\end{equation}

\noindent where $c_{27}:=c_{12} \, \max\{c_{25},c_{26}\}.$

 Recalling that $\sum_{p=1}^{M}
\gamma_{p}^{-\varepsilon'/6}\le c_{28} \,
\gamma_{M}^{-\varepsilon'/6}$ and $\sum_{p=1}^{M}
  \gamma_{p}\le c_{29} \,  \gamma_{1}=o(j_k)$ for some $c_{28},c_{29}>0,$ (\ref{rwrers:mino4}) yields
\begin{equation}
    \label{rwrers:mino5}
    P_z \left\{G(j_k, \ell) \, , \, H(j_k,\ell) \right\} \ge c_{30} \, Q_0  \,  P_z \left\{
    H(j_k,\ell) \right\},
\end{equation}

\noindent for some $c_{30}>0.$ Observe that Donsker's theorem
implies that there exists $c_{31}>0$ such that $\min\{P_z \{
H(j_k,\ell)\},Q_0\} \ge c_{31}.$ Therefore, assembling
(\ref{rwrers:mino3}) and (\ref{rwrers:mino5}), we get
$$
 P^+_{k,\ell} \ge c_{32} \, {\gamma_{M} \over \gamma_0},
$$

\noindent where $c_{32}:=c_{24} \, c_{30} \, c_{31}^2.$

Recalling (\ref{rwrers:PF+1}) and $\theta_k=\lfloor \, j_k /
\gamma_{M} \rfloor -1,$ we get, uniformly in $z \in I_{j_{k-1}},$
$$
P_z \left\{ E^{env}(j_k) \right\} \ge c_{32} \, \theta_k \,
{\gamma_{M} \over \gamma_0} \ge c_{33},
$$

\noindent for all large $k$ and for some $c_{33}>0,$ which concludes
the proof of (\ref{rwrers:F+}).
 \hfill$\Box$

\bigskip
\bigskip
\noindent {\Large\bf Acknowledgements}

\bigskip

I would like to thank my Ph.D. supervisor Zhan Shi for helpful
discussions. Many thanks are due to two anonymous referees for
careful reading of the original manuscript and for invaluable
suggestions.

\bigskip

%%
%% Authors's adresses
%%

{\footnotesize

\baselineskip=12pt

\hskip110pt Laboratoire de Probabilit\'es et Mod\`eles Al\'eatoires

\hskip110pt Universit\'e Paris VI

\hskip110pt 4 place Jussieu

\hskip110pt F-75252 Paris Cedex 05

\hskip110pt France

\hskip110pt {\tt zindy@ccr.jussieu.fr}

}%end of authors' adresses

\end{document}